\documentclass[]{article}

     \usepackage[final]{neurips_2020}

\usepackage[utf8]{inputenc} 
\usepackage[T1]{fontenc}    
\usepackage{hyperref}       
\usepackage{url}            
\usepackage{booktabs}       
\usepackage{amsfonts}       
\usepackage{nicefrac}       
\usepackage{microtype}      

\usepackage[toc,page]{appendix}
\usepackage{minitoc}

\usepackage{tocloft}

\usepackage{amsmath}
\usepackage{amsfonts,amssymb,amsthm}
\usepackage{graphicx}
\usepackage[dvipsnames]{xcolor}

\usepackage{mathrsfs,booktabs}
\usepackage{nicefrac}

\definecolor{darkmidnightblue}{rgb}{0.0, 0.2, 0.4}
\definecolor{darkpowderblue}{rgb}{0.0, 0.2, 0.6}
\definecolor{dukeblue}{rgb}{0.0, 0.0, 0.61}

\hypersetup{
    colorlinks = true,
    citecolor=blue,
    urlcolor=blue,
    breaklinks=true,
    linkcolor = dukeblue,
    linkbordercolor = {white},
}
\usepackage{natbib}
\usepackage{cleveref}
\usepackage{autonum}

\newcommand\bs{\boldsymbol}

\newcommand\bfE{\mathbf E}

\newcommand\bL{\boldsymbol L}
\newcommand\bY{\boldsymbol Y}
\newcommand\bP{\boldsymbol P}
\newcommand\bQ{\boldsymbol Q}

\newcommand\bX{\boldsymbol X}

\newcommand\bW{\boldsymbol W}

\newcommand\bx{\boldsymbol x}
\newcommand\by{\boldsymbol y}

\newcommand\NN{\mathbb N}
\newcommand\RR{\mathbb R}
\newcommand\EE{\mathbb E}

\newcommand\bfI{\mathbf I}
\newcommand\bfA{\mathbf A}
\newcommand\bfB{\mathbf B}

\newcommand\btheta{\boldsymbol\theta}
\newcommand\bvartheta{\boldsymbol\vartheta}

\def\tilde{\widetilde}

\def\bfE{\mathbf E}

\newcommand{\SCGL}[2]{\hyperlink{SCGL}
{{\textupsf{{$({#1},{#2})$-SCGL}}}}}

\newcommand{\textupsf}[1]{\textup{\textsf{#1}}}
\def\CLPEN1{\hyperlink{PEN}{\textupsf{{$(C,\lambda)$-PEN}}}}

\def\PLD1{\hyperlink{PLD}{\textupsf{{PLD}}}}
\def\PGF1{\hyperlink{PGF}{\textupsf{{PGF}}}}
\def\TPLD1{\hyperlink{TPLD}{\textupsf{{TPLD}}}}
\newcommand{\AssA}[2]{\hyperlink{AssA}{\textup{\textbf{A}}}$({#1},{#2})$}

\newtheorem{theorem}{Theorem}
\newtheorem{proposition}{Proposition}

\newtheorem{lemma}{Lemma}

\crefname{cor}{Corollary}{Corollaries}
\crefname{condition}{Condition}{C}

\date{September 2019}

\begin{document}

\doparttoc 
\faketableofcontents 

\part{} 

\title{Penalized Langevin dynamics with vanishing penalty 
for smooth and log-concave targets}

\author{%
  Avetik Karagulyan \\
  CREST, ENSAE, IP Paris\\
  \texttt{avetik.karagulyan@ensae.fr} \\
  \And
  Arnak S.\ Dalalyan \\
  CREST, ENSAE, IP Paris\\
  \texttt{arnak.dalalyan@ensae.fr}
}

\maketitle

\begin{abstract}
We study the problem of sampling from a probability distribution
on $\mathbb R^p$ defined via a convex and smooth potential function.  
We first consider a continuous-time diffusion-type process, termed 
Penalized Langevin dynamics (PLD), the drift of which is the negative 
gradient of the potential plus a linear penalty that vanishes when 
time goes to infinity. An upper bound on the Wasserstein-2 distance 
between the distribution of the PLD at time $t$ and the target is
established. This upper bound highlights the influence of the speed 
of decay of the penalty on the accuracy of approximation.
As a consequence, considering the low-temperature
limit we infer a new nonasymptotic guarantee of convergence of 
the penalized gradient flow for the optimization problem.
\end{abstract}

\setlength{\parindent}{0pt}

\section{Introduction}

The problem of sampling from a probability distribution received a great deal of attention in machine learning literature. Gradient based MCMC methods such as the Langevin MC, the underdamped Langevin Monte Carlo, the Hamiltonian Monte Carlo and their Metropolis adjusted 
counterparts were shown to have attractive features both in practice
and in theory. In particular, thanks to a large number of recent results, the case of smooth and strongly log-concave densities is now fairly well understood. In this case, non-asymptotic theoretical
guarantees for various distances on probability distributions have
been established, showing that the number of gradient evaluations
necessary to achieve an error upper bounded by $\varepsilon$ 
is a low order polynomial of the dimension, the condition number
and the inverse precision $1/\varepsilon$. The dependence on the latter is even logarithmic for Metropolis adjusted methods. 

The main focus of this paper is on the problem of sampling from 
densities\footnote{We will use the same notation for the probability density functions and corresponding distributions.}  
\begin{equation}
	\pi(\btheta) \propto \exp(-f(\btheta)),\qquad \btheta\in\mathbb R^p,
\end{equation}
corresponding to a (weakly) convex potential function $f:\mathbb R^p \to\mathbb R$. In the sequel, a twice differentiable function 
$f:\RR^p \rightarrow \RR$ is said to satisfy \SCGL{m}{M} condition \hypertarget{SCGL1}{($m$-strongly convex and $M$-gradient Lipschitz)}, for some $M \ge m\ge 0$, if the following inequality is satisfied:
\begin{equation}
		m \bfI_p \preceq \nabla^2 f(\btheta) \preceq M \bfI_p,\qquad \forall \btheta\in\mathbb{R}^p.
\end{equation}
Here, for two square matrices $\bfA$ and $\bfB$, the relation 
$\bfA\preceq \bfB$ means that $\bfB-\bfA$ is positive semidefinite. 

In this work, we wish to define a class of continuous-time
processes, such that for every element $\{\bL_t:t\ge 0\}$ of the class, 
the distribution of $\bL$ at time $t$ is close to the target 
distribution $\pi$. When the potential function $f$ satisfies
\SCGL{m}{M} condition with $m\ge 0$, it is well-known that the vanilla Langevin dynamics $\bL^{\textsf{LD}}$, defined as the solution of 
\begin{equation}\label{eq-langevin}\tag{\textsf{LD}}
	d \bL_t^{\textsf{LD}} = - \nabla f(\bL_t^{\textsf{LD}})dt 
	+ \sqrt{2}\,d\bs W_t, 
\end{equation}
where $\bs W_t$ is a standard Wiener process independent of $\bL^{\textsf{LD}}_0$, has $\pi$ as invariant distribution 
\citep{bhattacharya1978}. Furthermore, when $m>0$, the distribution
$\nu_t^{\textsf{LD}}$ of $\bL_t^{\textsf{LD}}$ converges in Wasserstein 
distance (see below for a definition) exponentially fast to $\pi$ 
\citep{villani2008optimal}, 
that is
\begin{equation}\label{eq:expLD}
 	W_2(\nu_t^{\textsf{LD}},\pi) 
 	    \leq e^{-mt} W_2(\nu_0^{\textsf{LD}},\pi).
\end{equation}
A remarkable feature of this result is that it is dimension free. 
In the case $m=0$, it was established by \citep{Bobkov99} that 
the target distribution satisfies the Poincar\'e inequality 
with the Poincar\'e constant $\mathcal C_\textsf{P}$ that might depend  on the dimension. According to \citep{chewi2020exponential}, 
the exponentially fast convergence to zero holds true with $m$ 
replaced by $1/{\mathcal C_\textsf{P}}$, when $m=0$. In \citep{Kannan95}, the authors conjectured that there is a universal
constant $\mathcal C_{\textsf{KLS}}>0$ such that for any log-concave distribution $\pi$ 
on $\mathbb R^p$, $\mathcal C_{\textup{P}} \le  
\mathcal C_{\textsf{KLS}} \|\mathbf E_{\bX 
\sim \pi}[\bX\bX^\top]\|_{\textup{op}}$, where $\|\bfA\|_{\rm op}$ stands for the operator norm of the matrix $\bfA$. Despite important efforts 
made in recent years (see \citep{KLS_conjecture, cattiaux2018poincar}),
this conjecture is still unproved. Note also that the Poincar\'e 
constant is, in general, hard to approximate and to estimate \citep{pillaudvivien2019statistical}.  

One approach to getting more tractable convergence bounds could be 
to find a tractable upper bound on the Poincar\'e constant of a 
distribution defined by a potential satisfying \SCGL{m}{M} condition
with $m=0$. We develop here another approach, consisting in modifying
the Langevin dynamics, so that the new dynamics has still a limiting
distribution equal to $\pi$ but for which we can get a tractable
upper bound. A natural way of defining this new dynamics is to add
to the potential $f$ a strongly-convex penalty with a strong-convexity
constant that vanishes when time goes to infinity. For a quadratic
penalty function, this leads to the process $\bL^{\textsf{PLD}}$
termed \hypertarget{PLD}{penalized Langevin dynamics} and defined by\footnote{If a 
good initial guess $\btheta_0$ of the minimum point of $f$ is
available, it might be better to replace the penalty term by
$\alpha(t)(\cdot - \btheta_0)$.}
\begin{equation}\label{eq-plangevin}\tag{\textsf{PLD}}
	d \bL_t^{\textsf{PLD}} = - \big(\nabla f(\bL_t^{\textsf{PLD}}) 
	+ \alpha(t) \bL_t^{\textsf{PLD}}\big)\,dt 
	+ \sqrt{2}\,d\bs W_t, 
\end{equation}
where $\alpha:[0,\infty)\to[0,\infty)$ is a time-dependent penalty 
factor tending to zero as $t\to\infty$. The main result of this work is an upper bound on 
$W_2(\nu_t^{\textsf{PLD}},\pi)$
that is valid for every continuously differentiable and decreasing penalty factor $\alpha$. Optimizing over $\alpha$, we show that the choice $\alpha(t) \sim 1/(2t)$, 
when $t\to\infty$, leads to a simple upper bound of the order $1/\sqrt{t}$. 
Interestingly, using a suitably parametrized temperature-dependent potential function 
$f_\tau(\cdot) = f(\cdot)/\tau$ and a penalty factor, one can get an upper bound for 
the penalized gradient flow
\begin{align}\tag{\hypertarget{PGF}{\textsf{PGF}}}
    \dot{\bX}{}_t^{\textsf{PGF}} = -\big(\nabla f(\bX_t^{\textsf{PGF}}) + \alpha_0(t)
    \bX_t^{\textsf{PGF}} \big),\qquad t\ge 0,
\end{align}
by  passing to the low-temperature limit.
This bound implies that $\| \bX_t^{\textsf{PGF}} - \bx^* \|_2$ can be 
of the order $O(1/t^{1-\mathsf q})$, where $\bx^*$ is a minimizer 
of the potential $f$ and $\mathsf q\in [0,1]$ is a parameter appearing
in an additional condition imposed on $f$. To the best of our knowledge,
the obtained bound is new, most previous results being valid for the
objective function itself, not for the minimum point.

The rest of the paper is organized as follows. 
We start by stating the bound on the error of the 
\PLD1 in \Cref{sec:PLD}. We also instantiate
the bound to the penalty factors that are 
inversely proportional to time. In \Cref{sec:PGF}, we 
discuss the connections with the optimization problem, 
assessing the error of the \PGF1. 
\Cref{sec:prior} is devoted to relation to 
prior work. 
The proof of the main result, up to some 
technical lemmas, is presented in \Cref{sec:proof}. Missing proofs are deferred 
to the supplementary material. 

To complete this introduction, we introduce some
notations. We consider   the Wasserstein-2 
distance
\begin{equation}
	 W_2(\nu,\nu')  = \inf\Big\{\bfE[\|\bvartheta-
    \bvartheta'\|_2^2]^{1/2}
    :\bvartheta\sim \nu \text{ and }
    \bvartheta'\sim \nu' \Big\},
\end{equation}
where the minimum is over all joint distributions
having $\nu$ and $\nu'$ as the first and the second marginal distributions. For any $\gamma 
> 0$, we define the probability density function $\pi_\gamma (\btheta) \propto \exp(-f(\btheta) - \gamma \|\btheta\|_2^2)$, where $\|\btheta\|_2$ 
is the Euclidean norm. We also define $\mu_k(\pi) = \EE_{\bX\sim \pi}[\|\bX\|_2^k]$, the moment
of order $k$ of $\pi$. 

\section{Convergence of penalized Langevin dynamics}\label{sec:PLD}

In this section we explain our approach and state the main result. 
Without loss of generality, we will assume that the initial point
for the \PLD1 is the origin, $\bL_0^{\textsf{PLD}} = 0$. Note that if
a good guess $\btheta_0$ of a minimizer of $f$ is available, it is
recommended to initialize \PLD1 at $\btheta_0$. Our framework covers
this case, since it suffices to apply our results to the translated
function $\tilde f(\cdot) = f(\btheta_0 + \cdot)$. Under the condition
$\bL_0^{\textsf{PLD}} = 0$, the Wasserstein-2 distance at the starting
point coincides with the second-order moment, $W_2(\nu_0^{\textsf{PLD}},\pi) = 
\sqrt{\mu_2(\pi)}$.

When $\alpha(t) = \alpha$ is a strictly positive constant, the
distribution $\nu_T^{\textsf{PLD}}$, for a large value of $T$, is 
close to the biased target $\pi_\alpha$. Furthermore, in view of
\eqref{eq:expLD}, the distance between these two distributions is 
smaller than a prescribed error level $\varepsilon>0$ as soon as 
$T \ge (1/\alpha)\log(\sqrt{\mu_2(\pi)}/\varepsilon)$. On the other 
hand, one can choose $\alpha$ small enough such that the bias $W_2(\pi_\alpha,\pi)$ is smaller than $\varepsilon$. The discrete 
counterpart of this approach has been used in many recent works \citep{Dalalyan14,dalalyan2019bounding,Dwivedi}.  The approach we
develop here extends these work to the case of time-dependent $\alpha$ 
and has the advantage of being asymptotically unbiased, when 
$t\to\infty$. In other words, it allows to choose $\alpha$ 
independently of $\varepsilon$ and make the error smaller than
$\varepsilon$ by running \PLD1 over a sufficiently large time period. 

\begin{theorem}\label{th-main} 
	Suppose that $\pi$ is a probability distribution with 
	a potential function $f$ that satisfies 
	\SCGL{m}{M} condition, where
	$m\geq 0$  and $M>0$. Let $\alpha : [0,+\infty) \rightarrow \RR$ 
	be a non-increasing differentiable function, such that $  m +
	\alpha(t) >0$ for every $t\ge 0$. Then, for every positive number 
	$t$ and for $\beta(t) = \int_0^t (m+ \alpha(u))\, du$,  we have
    \begin{align}\label{main:ineq}
	    W_2(\nu_t^{\textupsf{{PLD}}}\!,\pi)
		& \leq  \sqrt{\mu_2(\pi)}\, e^{-\beta(t)}  + 11\mu_2(\pi)\bigg\{
		e^{\beta(t)}\int_0^t \frac{|\alpha'(s)|e^{\beta(s)}}{\sqrt{m+
		\alpha(s)}}\, ds	+	\frac{ \alpha(t)}{\sqrt{m+\alpha(t)}}\bigg\}.
    \end{align}
\end{theorem}

The proof being postponed to \Cref{sec:proof}, the rest of this section
is devoted to discussing the stated theorem and its consequences for
some specific choices of the penalty factor. One can notice right away 
that in the case of a positive $m$, we can choose $\alpha$ to be zero, thereby obtaining the classical exponential convergence rate 
\citep{villani2008optimal}. In the rest of the discussion, we assume 
that $m=0$.

The numerical constant 11 can certainly be improved. It is closely
related to the fact that for a log-concave distribution $\nu$, we 
have $\mu_4(\nu) \le C_4 \mu_2^2(\nu)$ for a universal constant $C_4$. 
proved to satisfy $C_4\le 442$ \citep[Remark 3]{dalalyan2019bounding}. 
Improved bounds on $C_4$ will automatically yield improved numerical
constant in \Cref{th-main}. We also note that the Lipschitz constant
$M$ does not appear in inequality \eqref{main:ineq}. Our proof, however,
 requires the finiteness of $M$. 
We believe that it is possible to relax the gradient-Lipschitz assumption
by requiring from $\nabla f$ to be only locally Lipschitz-continuous. 

One can check that if we replace the penalty factor $\alpha$ by a larger
function, the first term of the upper bound in \eqref{main:ineq}, proportional to $\exp\{-\int_0^t \alpha(s)\,ds\}$, becomes smaller. 
On the other hand, the second term increases when $\alpha$ increases
\footnote{This is clear for the last term, proportional to $\sqrt{\alpha(t)}$, whereas the corresponding claim for the first term
in the curly parentheses less trivial.} One can thus use inequality 
\eqref{main:ineq} for choosing the penalty factor $\alpha$ that offers
a trade-off between the error of approximating the biased density
$\pi_{\alpha(t)}$ by the \PLD1 and the error of approximating the target
$\pi$ by the biased density $\pi_{\alpha(t)}$. 

To ``optimize'' the upper bound with respect to $\alpha$, 
let us for the moment ignore the first term in the curly parentheses in 
\eqref{main:ineq}. In that case our functional of interest has two
components, where one of them is increasing with respect to $\alpha$,
while the other is decreasing. (Here the monotony must be understood with a  certain precaution, as in our case the mathematical concept is
not well-defined.) These considerations suggest to choose the ``optimal" 
$\alpha$ by balancing these two terms:
\begin{equation}\label{eq-first-third}
 	\sqrt{\mu_2(\pi)}\, e^{-\beta(t)} = 11\sqrt{\alpha(t)}\,  \mu_2(\pi).
\end{equation}
Taking the square of both sides, cancelling out some terms and using
that $\alpha(t) = \beta'(t)$, we check that \eqref{eq-first-third} is equivalent to 
\begin{align}
	121\mu_2(\pi)\beta'(t)e^{2\beta(t)} = 1.
\end{align}
Solving this differential equation we get the following expression 
for $\beta(\cdot)$ and the corresponding expression for $\alpha(\cdot)$:
\begin{equation}
\beta^*(t) = \frac{1}{2} \log\left(\frac{2t}{121\mu_2(\pi)}+1\right) \hspace{0.4cm} \text{and} 
\hspace{0.4cm} \alpha^*(t) = \frac{1}{2t + 121\mu_2(\pi)}.
\end{equation} 
It is easy to check that this choice of $\alpha$ ensures that the
first and the last terms in the right hand side of \eqref{main:ineq}
are of the order $1/\sqrt{t}$. Interestingly, the middle term in
\eqref{main:ineq} turns out to be of the same order, up to a logarithmic
factor. The precise statement of the consequence of \Cref{th-main} 
obtained by choosing $\alpha(t) = 1/(2t+A)$ for some $A> 0$ reads as
follows.
\begin{figure}
    \centering
    \hspace{-13pt}
    \includegraphics[height = 126pt]{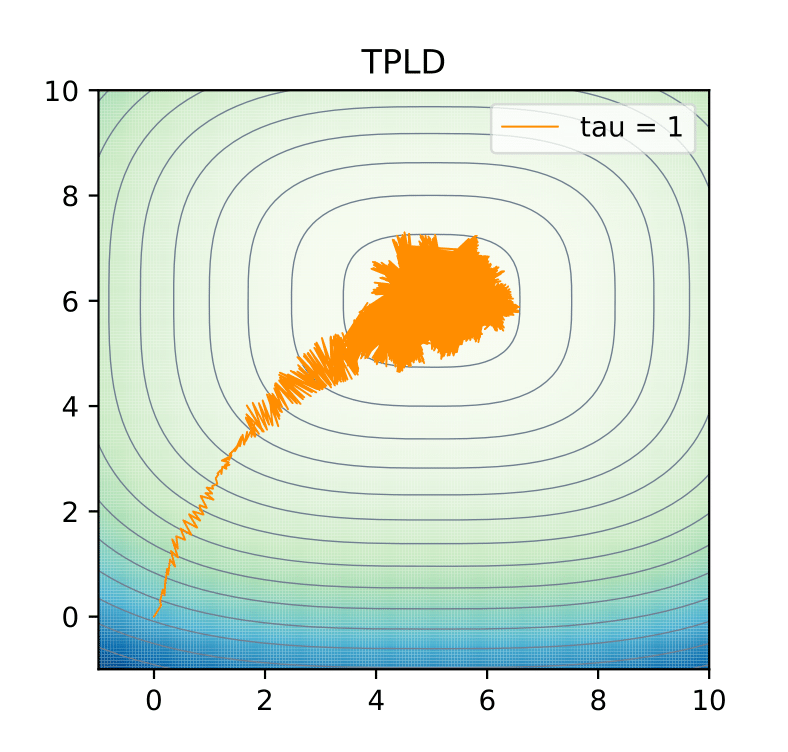}\hspace{-10pt}
    \includegraphics[height = 126pt]{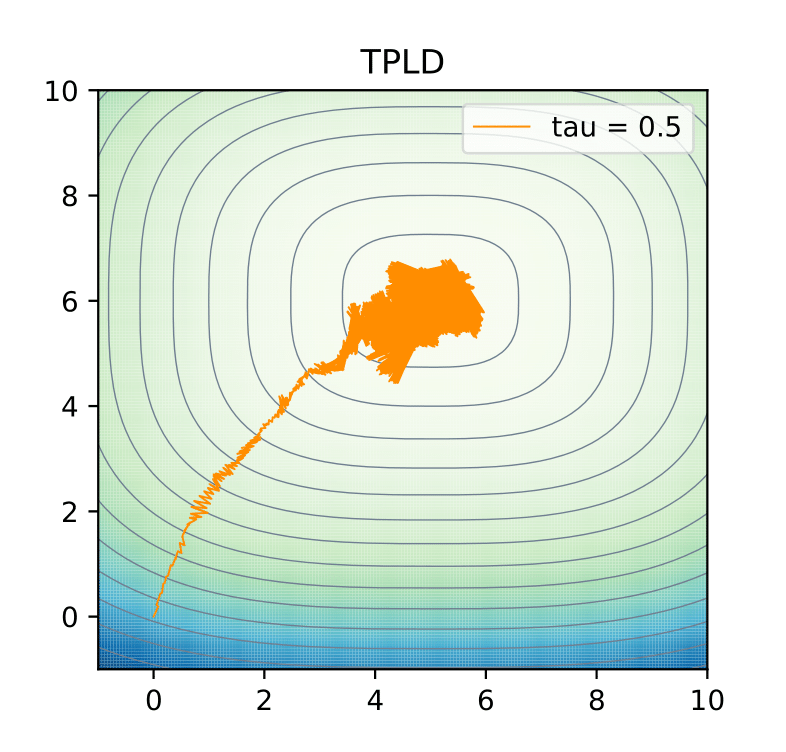}\hspace{-15pt}
    \includegraphics[height = 126pt]{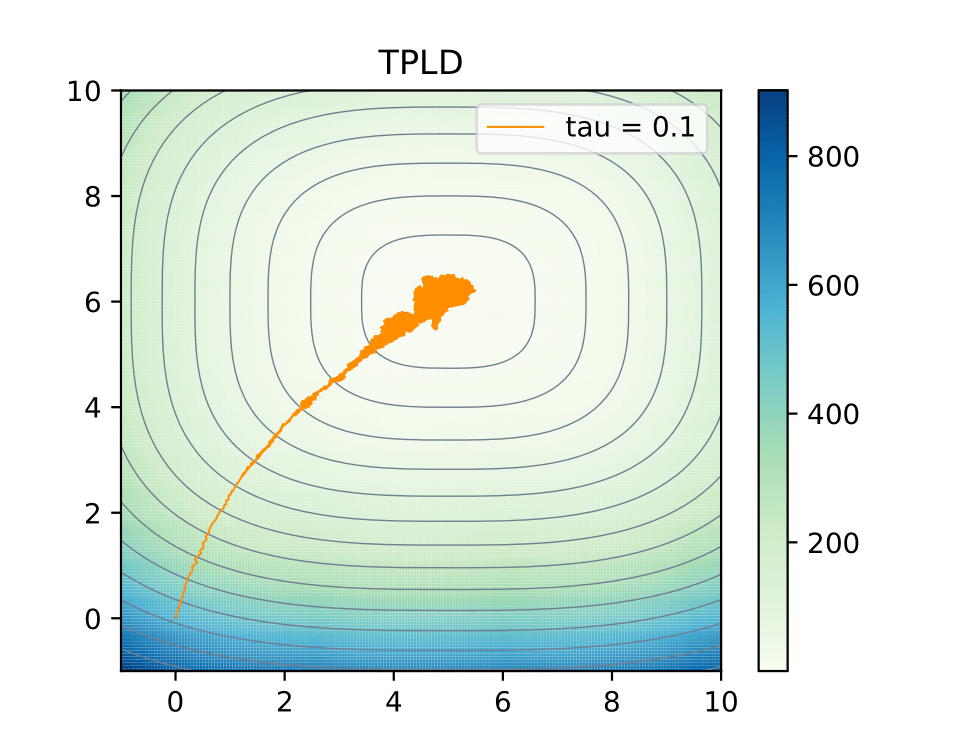}
    \vspace{-5pt}
    \caption{One random path of the \protect\TPLD1, for $\tau=1$ (left), 
    $\tau = 0.5$ (middle) and $\tau =0.1$ (right). }
    \vspace{-5pt}
    \label{fig:0}
\end{figure}
\begin{proposition}\label{prop-optimal}
	If the potential function satisfies the \SCGL{m}{M} condition
	with $m=0$, then the error of \PLD1 with $\alpha(t) = 1/({A + 2t})$, measured by the Wasserstein-2 distance, satisfies
	\begin{equation}
		 W_2(\nu_t^{\textupsf{{PLD}}} , \pi)   \leq \frac{\sqrt{A\mu_2(\pi)} + 11 \mu_2(\pi)\big(1 + 
		 \log{\left(1 + (2/A)t\right)}\big)}{\sqrt{A + 2t}}.
	\end{equation}
	In particular, for $A = 2\mu_2(\pi)$, we get
	\begin{equation}
	W_2(\nu_t^{\textupsf{{PLD}}} , \pi)  \leq
	    \frac{10\mu_2(\pi)\big\{1 + \log{\big(1 + (t/\mu_2(\pi))\big)} 
	    \big\}}{\sqrt{\mu_2(\pi) + t}}.
	\end{equation}
\end{proposition}

To complete this section, we present a quick argument showing that
the right hand side of \eqref{main:ineq} cannot converge to zero faster 
than the rate $1/\sqrt{t}$. Suppose that for a specific choice of 
$\alpha$, the right hand side of \eqref{main:ineq} is $o(1/t^{1/2})$.
Then the last term is necessarily $o(1/t^{1/2})$, which implies that
$\alpha(t) = o(1/t)$  when $t\to +\infty$. Thus, for some $c > 0$,
$\alpha(t)\le 1/(4t)$ for every $t\ge c$. This means that $\beta(t) 
\le c\alpha(0) +  (\nicefrac14) \log(t/c)$. Hence,
\begin{equation}
	\exp(-\beta(t)) \geq  {(c/t)^{1/4}}\,{\exp(- c\alpha(0))}.
\end{equation}
This proves that the upper bound on the error of the 
\PLD1 established in \Cref{th-main} may not tend to zero at a 
faster rate than $1/\sqrt{t}$, as $t$ goes to infinity. 
This argument also shows that the optimizer $\alpha(\cdot)$ of the upper 
bound is asymptotically equivalent to $1/(2t)$ when $t\to\infty$.

\section{The counterpart in optimization: penalized gradient flow}\label{sec:PGF}

In this section we draw the parallel between \PLD1 and the penalized gradient flows, 
henceforth referred to as \PGF1, for a non-strongly convex function $f$. In the case of strongly convex functions, the gradient flow $\bX^{\sf GF}_t$ defined by the 
differential equation $\dot{\bX}{}^{\sf GF}_t = -\nabla f(\bX_t^{\sf GF})$, converges 
exponentially fast to the minimum $\bx_*$ of $f$, without the need 
to add a quadratic penalty. In contrast with this, for general non-strongly 
convex case functions $f$, only the convergence of the function $f(\bX_t^{\sf GF})$ to
$f(\bx_*)$ at the rate $1/t$ can be established. The goal
of this sections is to understand the convergence of the flow to the minimum point 
when a vanishing
quadratic penalty is added to the cost function $f$; when does this flow converge, 
what is the impact of the penalty factor and what kind of
rate can be achieved. To answers these questions, we assume in this section that 
\SCGL{0}{M} holds true. We also assume that $f$ has a 
unique minimum point denoted by $\bx_*$. 

In the analysis performed in the previous section, we can replace the function 
$f(\cdot)$ by the function $f_{\tau}(\cdot)= f(\cdot) / 
\tau$. The function $f_{\tau}$ has $\bx_{*}$ as its point of minimum,  whatever 
the real number $\tau>0$. Moreover, if we define the tempered density function 
$\pi^{\tau}({\btheta})
\propto \exp \left(-f_{\tau}(\btheta)\right)$,
the distribution $\pi^{\tau}$ tends to  $\delta_ 
{\bx_*}$, the Dirac mass at $\bx_*$. Clearly,
$f_{\tau}$ satisfies \SCGL{0}{M/\tau} condition. 
Thus, according to \Cref{th-main}, the process $\bL^\tau_t$, defined as 
\begin{align}\label{eq-L-tau}
    \bL^\tau_t = \bL^\tau_0 - \frac{1}{\tau}\int_0^{t}\big( \nabla 
    f(\bL^\tau_s) + \alpha(s/\tau)\bL^\tau_s\big) ds + \sqrt{2}\,\bW_t,
\end{align}
converges to $\pi^\tau$ in Wasserstein distance, 
if  $\alpha(\cdot)$ is a continuously differentiable and non-increasing function. 
We now introduce the \hypertarget{TPLD}{tempered penalized Langevin dynamics 
(\textsf{TPLD})},  as a time-scaled version of $\bL^\tau$:   $\bX^{\textsf{TPLD}}_t = 
\bL^\tau_{t\tau}$ for every $\tau > 0$. One can check that this process satisfies the 
stochastic differential equation 
\begin{align}\label{eq:TPLD}\tag{\textsf{TPLD}}
    d\bX^{\textsf{TPLD}}_t = -  \big( \nabla 
    f(\bX^{\textsf{TPLD}}_t) + \alpha(t) \bX^{\textsf{TPLD}}_{t}\big)\,dt + \sqrt{2\tau}\, d\bar{\bW}_{t},
\end{align}
where $\bar\bW_t = \tau^{-1/2}\bW_{\tau t}$ 
is a standard Wiener process. To illustrate the behaviour of this process,
\Cref{fig:0} shows one realization of \TPLD1 for different values of $\tau$, 
with the left plot corresponding to \PLD1. All the results of the previous
section continue to hold for this tempered dynamics. In particular, the
analog of the second claim of \Cref{prop-optimal} in the case of the tempered 
diffusion takes the following form.

\begin{proposition}\label{prop-optimal-tempered}
	If the potential function satisfies the \SCGL{0}{M} condition, 
	then the error of\/ \TPLD1 with $\alpha(t) = 1/({2\mu_2(\pi^\tau) 
	+ 2t})$ satisfies
	\begin{equation}
	W_2(\nu_t^{\textupsf{{TPLD}}} , \pi^\tau)  \leq
	    \frac{10\mu_2(\pi^\tau)\big\{1 + \log{\big(1 + (t/\mu_2(\pi^\tau))\big)} 
	    \big\}}{\sqrt{\tau(\mu_2(\pi^\tau) + t)}},\qquad\forall t\ge 0.
	\end{equation}
\end{proposition}

As mentioned above, for small $\tau$, $\pi^\tau$ is close to the
Dirac mass at the minimum point $\bx_*$. The last result tells us
that we can approximate $\pi^\tau$ arbitrarily well, by running
the \TPLD1 over a large time-period. But we can not replace $\tau$ 
by zero in this result, since the denominator of the right hand 
side vanishes and the result becomes vacuous. We show below that this
can be repaired if an additional assumption is introduced. 

Taking $\tau = 0$ in \eqref{eq:TPLD}, we get
the \hypertarget{PGF}{penalized gradient flow (PGF)}:
\begin{align}\label{eq-pflow}
    d\bX_t^{\textsf{PGF}} = -\big( \nabla 
    f(\bX_{u}^{\textsf{PGF}}) + \alpha(u) \bX_{u}^{\textsf{PGF}}\big)\, du,\qquad t\ge 0,\qquad \bX_0^{\textsf{PGF}} = \mathbf 0.
\end{align}
Here we recognize the analog of \PLD1 in the setting of the gradient
flows. On the other hand, the Euclidean distance on $\RR^p$ is equal 
to the Wasserstein distance between Dirac measures. This leads us to
think that our approach for obtaining non-asymptotic error bounds for
\PLD1 is applicable to the penalized gradient flow. This turns out to
be true, modulo the introduction of the following assumption. 

\hypertarget{AssA}{\textbf{Assumption A$(\mathsf D, \mathsf q)$}:}
The minimum point $\bx_\gamma $ of the (strongly convex and coercive)
function $f_\gamma(\cdot) = f(\cdot) + \gamma\| \cdot \|_2^2/2$
satisfies 
\begin{align}\label{AssA}
    \|\bx_{\gamma} - \bx_{\tilde\gamma}\|_2 \leq  \frac{\mathsf D}{\tilde\gamma^{\mathsf q}}(\tilde\gamma-\gamma)
    \|\bx_*\|_2,\qquad \forall \gamma < \tilde\gamma,
\end{align}
for some $\mathsf D>0$ and $\mathsf q \in [0,1]$.  

Since $\bx_\gamma$ a stationary point of $f_\gamma$, we have 
$\nabla f(\bx_\gamma) = -\gamma \bx_\gamma$. From this relation and 
\citep[Theorem 2.1.12]{Nest}, one can deduce that (a) if $f$ satisfies
\SCGL{m}{M} condition with $m>0$ then \AssA{1/m}{0} holds  and 
(b) if $f$ satisfies \SCGL{0}{M} condition then \AssA{1}{1} 
holds.

\begin{theorem}\label{th-flow}
    Assume that $\alpha:[0,\infty)\to[0,\infty)$ is a continuously 
    differentiable and non-increasing function. Let $\beta(t) = \int_0^t \alpha(s)\,ds$ be the antiderivative of $\alpha$. If $f$ satisfies
    \AssA{\mathsf D}{\mathsf q} and \SCGL{0}{M}, then
    \begin{align}\label{bound:PGF}
        \|\bX_t^{\textupsf{{PGF}}} - \bx_*\|_2
    	& \leq   \|\bx_ *\|_2 \bigg(e^{-\beta(t)}  + \textupsf{D} 
        \int_0^t\frac{|\alpha'(s)|}{{\alpha^{\textupsf{q}}(s)}}
        e^{\beta(s) -\beta(t)} ds+ \textupsf{D} \alpha^{1-\textupsf{q}} 
        (t)\bigg).
    \end{align}
\end{theorem}
The proof can be found in \Cref{sec-C} of the supplementary material. 
When $\textsf{q}= 1/2$, this result is the optimization counterpart 
of the inequality shown in \Cref{th-main}. Once again, it is appealing 
to optimize the right hand side of \eqref{bound:PGF} in order to choose
the ``best'' penalty factor. Using arguments similar to those of
previous section, \textit{i.e.}, balancing the first and the last
terms on the right hand side of \eqref{bound:PGF},  we get 
that the ``optimal'' convergence of \PGF1 is achieved when
\begin{align}
    \alpha^*(t) = \frac{(1-{\textupsf q}) }{t + 
    \textupsf{D}^{1/(1-\textupsf{q})}(1-\textupsf{q})} \qquad \text{and} \qquad  
    \beta^*(t) =(1-\textupsf{q})\log\left(
    \frac{t}{\textupsf{D}^{1/(1-\textupsf{q})}(1-\textupsf{q})} + 1\right).
\end{align}
\begin{figure}
    \centering
    \hspace{-15pt}
    \includegraphics[height = 110pt]{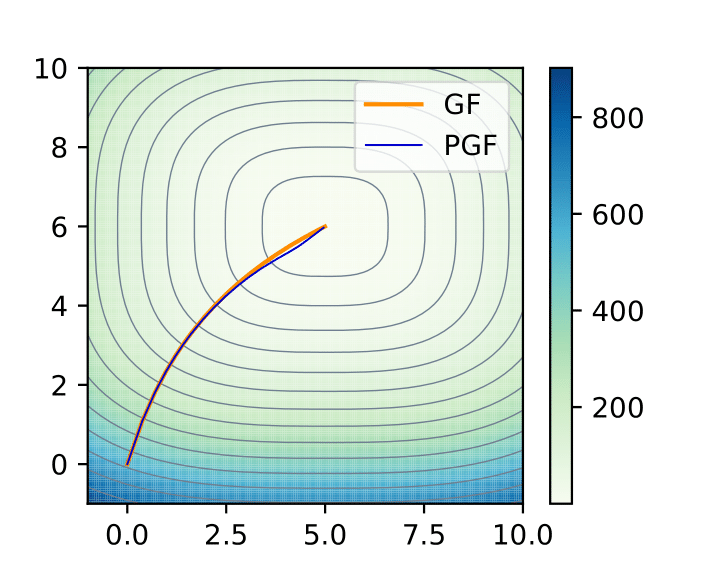}\hspace{-12pt}
    \includegraphics[height = 110pt]{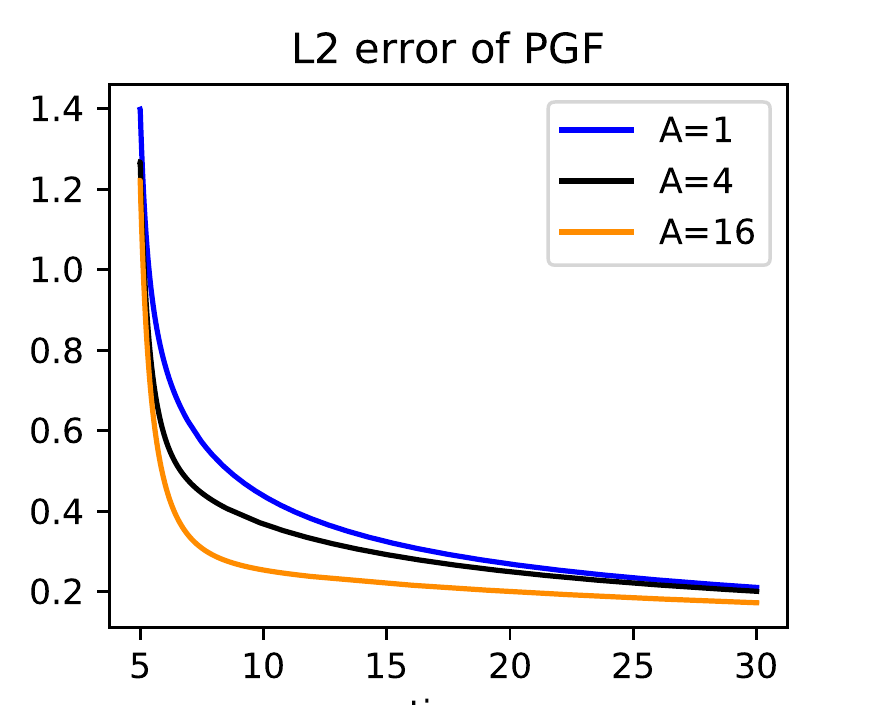}\hspace{-15pt}
    \includegraphics[height = 110pt]{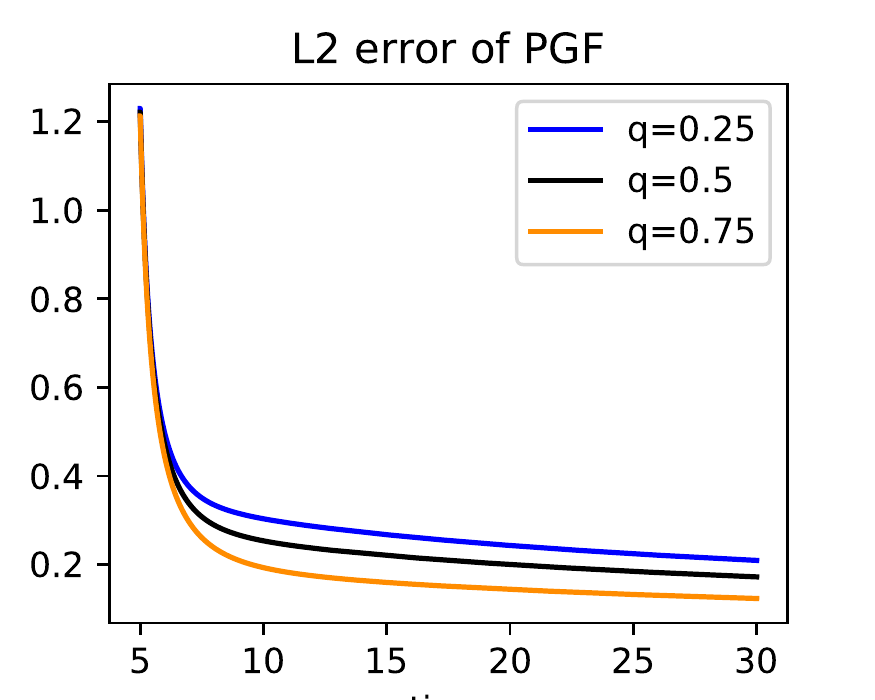}
    \caption{Left: The gradient flow (orange) and the penalized gradient flow (blue). Middle and Right: the mapping $t\mapsto \|\bX_t^{\textsf{PGF}}-\bx_*\|_2$. In both plots we used the function $f(x,y) = 
    |x-5|^3 + 2|y-6|^3$.
    \vspace{-10pt}}
    \label{fig:1}
\end{figure}
This leads us to make the recommendation of choosing 
$\alpha(t) = (1-\textsf{q})/(t+A)$, for some positive $A$. If $\mathsf q = 1$, this amounts to considering the non-penalized gradient flow
and \eqref{bound:PGF} boils down to the fact that the distance from the 
gradient flow of a convex function to its minimum decreases. While for $\mathsf q<1$,  for the foregoing choice of the penalty factor, we get the error bound 
\begin{equation}\label{bound:PGF2}
    \|\bX_t^{\textupsf{{PGF}}} - \bx_*\|_2\leq
    \frac{A^{1-\mathsf{q}} + \mathsf{D}+ \mathsf{D}\log(1+ (t/A))}{ 
    (t+A)^{1-\textupsf{q}}}\|\bx_*\|_2.
\end{equation}
To complete this section, we make some remarks on assumption \AssA{\mathsf D}{\mathsf q}. First, one can relax this assumption
by requiring that the desired inequality holds for sufficiently small
values of $\tilde\gamma$ only. In this form, it can be easily seen
that larger values of $\mathsf q$ correspond to weaker assumption. 
Second, even if the function $f$ is not strongly convex, it may
satisfy \AssA{\mathsf D}{0}. An example is the function 
$f(x) = \sqrt{(x-x_*)^2 + b^2}$. The second derivative of this function
is equal to $b^2/((x-x_*)^2+ b^2)^{3/2}$. This implies that $f$ 
satisfies \SCGL{0}{1/b}. We show in the supplementary material that 
it satisfies \AssA{\mathsf D}{0} for some finite value of 
$\mathsf D>0$. This is not really surprising, given that this
function is strongly convex on any compact set. Another instructive
example is the function $f(x) = |x-x_*|^{a}$ with $a \ge 2$. 
On compact sets, this function satisfies\footnote{See supplementary material.} \AssA{\mathsf D}{(a-2)/(a-1)}. Therefore, the error bound
\eqref{bound:PGF2} implies that \PGF1 converges to the minimum of $f$ 
at the rate $1/t^{1/(a-1)}$, which is faster than the rate derived
from the standard $O(1/t)$ bound for the non-penalized gradient flow. 
This behaviour is depicted in \Cref{fig:1}.



\section{Prior work and outlook}\label{sec:prior}

The relation of our results to some prior work has already been
highlighted in previous sections. This section provides some 
complementary bibliographic remarks on recent advances on 
Langevin diffusions, gradient flows and their discrete counterparts.

Convergence of Langevin dynamics in continuous time has received a lot of attention in probability, see \citep{Cattiaux1,Cattiaux2,Bolley} and
the references therein. An interesting known fact, for instance, is that
the Langevin dynamics satisfies\footnote{Here, $\nu_{t,\bx}^{\textsf{LD}}$
is the distribution at time $t$ of the \textsf{LD} starting at $\bx$, and the inequality is assumed to hold for any $\bx,\by\in\mathbb{R}^p$ and any $t\in\mathbb{R}$.} $W_2(\nu_{t,\bx}^{\textsf{LD}}, 
\nu_{t,\by}^{\textsf{LD}}) \le e^{-mt}\|\bx-\by\|_2$ if and only if
$f$ is $m$-strongly convex. More recently, many papers in statistics and machine learning literature established non-asymptotic error bounds for discretized algorithms, mainly focusing on the convex case, see \citep{Durmus2,durmus2017, Hsieh:266354, Bubeck18,shen2019randomized} 
in addition to previously cited papers. The non-convex
case was emphasized in \citep{Cheng3,majka2018nonasymptotic,Erdogdu18,
MangoubiV19}.

In the optimization setting, the results on the convergence of the
gradient flow for convex objectives have been known for a long time. 
More recently, \citep{SuBoydCandes} derived a continuous-time
second-order differential equation characterizing the Nesterov
acceleration. Continuous-time Accelerated Mirror Descent was 
studied in \citep{Krichene15}. An approach based on Bregman-Lagrangian 
functional for continuous-time momentum and other methods was 
developed in \citep{Wibisono,wilson2016lyapunov}. 
Further results on related topics, relevant to machine learning,
were obtained in \citep{Jingzhao18,Scieur,franca18a}. An overview of 
results on gradient flow beyond the Euclidean space setting can 
be found in \citep{Ambrosio,Santa17}. 

On a related note, several studies took advantage of
the fact that the distribution of the Langevin dynamics is a gradient 
flow in the space of measures \citep{Cheng1, Bernton18,
durmus2018analysis, wibisono18a}, in order to establish error bounds for
sampling algorithms. The relation with MMD was studied by \cite{Arbel19}. 

The results presented in the present work can be generalized in
various directions. In particular, it would be interesting to relax 
the smoothness assumption, following
an argument from, \textit{e.g.}, \citep{Durmus3,chatterji2019langevin,Mou,Salim}, 
to develop a similar theory for the kinetic Langevin dynamics \citep{eberle2019,Cheng2,dalalyan_riou_2018,YiAnMa} or to consider
the case of a Lévy process driven stochastic differential equation 
in the spirit of \citep{umut20,liang2019exponential}.

\section{Conclusion}

We put forward a family of time-inhomogeneous diffusion 
processes that converge to a pre-specified target distribution 
and, therefore, can be used for approximate sampling. These
processes are defined as penalized Langevin dynamics with a
penalty that vanishes when time goes to infinity. The penalty
allows to ensure strong convexity, which helps to handle
situations where the original log-density is not strongly 
convex. We established a simple non-asymptotic error bound
showing that the rate of convergence in the Wasserstein-2
distance is $o(1/\sqrt{t})$. We have also discussed analogous
results for the penalized  gradient flow. The important next 
step to investigate in future works is the analysis of 
discretized versions of the penalized Langevin dynamics.

\section{Proof of \Cref{th-main} 
}
\label{sec:proof}

Recall that for every $\gamma \in \RR $,  
$\pi_{\gamma}$ is the probability distribution with 
density proportional to 	$\exp(-f(\btheta) - \gamma \|\btheta\|^2_2/2)$. The  
triangle inequality for the Wasserstein distance yields
\begin{equation}\label{eq-first-tr-ineq}
W_2(\nu_t^{\textsf{PLD}},\pi) \leq 	W_2(\nu_t^{\textsf{PLD}},
\pi_{\alpha(t)}) + W_2(\pi_{\alpha(t)},\pi),
\end{equation} 
for every $t>0$.
We will bound these two terms separately, but let us start by 
stating	two technical lemmas. The first one is a consequence of 
the	well-known transportation cost inequality (see
\cite[Corollary 7.2]{gozlan2010transport}), whereas the second one
establishes the smoothness and the monotony with respect to $\gamma$ of the second-order moment of $\pi_\gamma$. The proofs
of these lemmas are postponed to \Cref{sec-appendix}.

\begin{lemma}\label{lem-tr-ineq} Let $\pi$ be a probability density function such that the potential $f = -\nabla \log \pi$ satisfies 
the \SCGL{m}{+\infty} condition.
Let $\tilde{\gamma} \geq \gamma$ be real numbers,  such that 
$m + \gamma\ge 0$. Then
\begin{align}
	 W_2(\pi_{{\tilde{\gamma}}}, \pi_\gamma)&\le  \frac{ 11({\tilde{\gamma}}- \gamma)}
	 {\sqrt{m+{\tilde{\gamma}}}}\mu_2(\pi_\gamma).\label{lem1:1}
\end{align}
\end{lemma}
 
\begin{lemma}\label{lem-mon-mu} 
	Suppose that $\pi$ has a finite fourth-order moment. Then 
	$\gamma \mapsto \mu_2(\pi_{\gamma})$ is 
	 continuously differentiable and  non-increasing, when 
	$\gamma \in  [0,+\infty)$. 
\end{lemma}

If we  apply \Cref{lem-tr-ineq} with $\gamma = 0$ and 
$\tilde{\gamma} = \alpha(t)$, then we obtain
\begin{equation}\label{eq-last-term}
	 W_2(\pi_{\alpha(t)} , \pi )\le \frac{ 11\alpha(t)}{\sqrt{m+\alpha(t)}}
\mu_2(\pi).
\end{equation}
This provides the desired upper bound on the second term of the right hand side of 
\eqref{eq-first-tr-ineq}.	To bound the first term, 
$W_2(\nu^{\textsf{PLD}}_t,\pi_{\alpha(t)})$, we aim at obtaining a
Gronwall-type inequality for the function 
\begin{equation}
    \phi(t):= W_2(\nu_{t}^{\textsf{PLD}},\pi_{\alpha(t)}).
\end{equation} 
To this end, we consider an auxiliary 
stochastic process $\{\tilde{\bL}_{u}:u\ge t\}$, defined as 
a solution of the following stochastic differential equation
\begin{equation}\label{def-L-tilde}
    d\tilde{\bL}_{u} = -\big(\nabla f(\tilde{\bL}_{u}) + \alpha(t) \tilde{\bL}_{u}\big)du + \sqrt{2} d\bs W_u,
\end{equation}
with the starting point $\tilde{\bL}_{t} = {\bL}_{t}^{\textsf{PLD}}$.  
This is in fact the Langevin diffusion corresponding to the potential $f(\cdot) + \alpha(t)^2\|\cdot\|_2^2/2$. Therefore, $\pi_{\alpha(t)}$ is the
invariant distribution of $\tilde{\bL}$, and it is 
$\left(m+\alpha(t)\right)$-strongly log-concave. 
Let $\bQ_{t,\delta}$ be the distribution of the random vector  
$\tilde{\bL}_{t+\delta}$. 	The triangle inequality 
yields
\begin{equation}
         \phi(t+\delta) 
         \leq 
         W_2\left(\nu_{t+\delta}^{\textsf{PLD}}, \bQ_{t,\delta}\right) +
         W_2\left( \bQ_{t,\delta}, \pi_{\alpha(t)} \right) + 
         W_2\left(\pi_{\alpha(t)}, \pi_{\alpha(t+\delta)} \right) .     
\end{equation}
Recalling the definition of $\pi_{\alpha(t)}$ and $\bQ_{t,\delta}$, we therefore find ourselves in 
the case of classical Langevin diffusion. Hence, one can apply  \eqref{eq:expLD} to get the bound
\begin{equation}\label{eq-contr-delta}
	W_2( \bQ_{t,\delta},\pi_{\alpha(t)} ) \leq 
		\exp\big(-\delta(m+ \alpha(t))\big)
	W_2( \nu_t^{\textsf{PLD}}, \pi_{\alpha(t)} ) 
		= \exp \left( -\delta ( m+ \alpha(t)) \right) \phi(t).
\end{equation}  
 Applying \Cref{lem-tr-ineq} to $\pi_{\alpha(t)}$ and $\pi_{\alpha(t+\delta)}$, we get
 \begin{equation}\label{eq-dist-t-t+delta}
 	 W_2\left(\pi_{\alpha(t)}, \pi_{\alpha(t+\delta)} \right) 
 	 \leq \frac{11(\alpha(t)- \alpha(t+\delta))}
   {\sqrt{m+\alpha(t)}}\mu_2(\pi_{\alpha(t+\delta)}).
 \end{equation}
Thus we obtain a bound  for $\phi(t+\delta)$, that depends linearly on $\phi(t)$:
\begin{align}\label{ineq-bef-division}
	\phi(t+\delta)\leq  
		 W_2\left(\nu_{t+\delta}^{\textsf{PLD}}, \bQ_{t,\delta}\right) 
	& +  e^{-\delta(m+ \alpha(t))}\phi(t) 
     + \frac{11(\alpha(t)- \alpha(t+\delta))}
   {\sqrt{m+\alpha(t)}}\mu_2(\pi_{\alpha(t+\delta)}).
\end{align}
Let us subtract $\phi(t)$ from both sides of 
\eqref{ineq-bef-division} and divide by $\delta$:
\begin{align}
	\frac{\phi(t+\delta) - \phi(t)}{\delta}\leq 
		\frac{1}{\delta}\cdot W_2\left(\nu_{t+\delta}^{\textsf{PLD}}, \bQ_{t,\delta}\right) 
	& + \frac{\exp\left(-\delta(m+ \alpha(t))\right) - 1}{\delta}\cdot \phi(t)\\
   & +\frac{11(\alpha(t)- \alpha(t+\delta))}
   {\delta\sqrt{m+\alpha(t)}}\mu_2(\pi_{\alpha(t+\delta)}).\label{eq-minus-phi}
\end{align}
The next lemma provides an upper bound on  
$W_2\left(\nu_{t+\delta}^{\textsf{PLD}}, \bQ_{t,\delta}\right) $ showing that it is $o(\delta)$, when $\delta \rightarrow 0$.
 \begin{lemma}\label{lem-gr-1st}
	For every  $t,\delta>0$, and for every integrable function $\alpha:[t,t+\delta]\to [0,\infty)$,
	\begin{equation}
      W_2\left( \nu_{t+\delta}^{\textupsf{{PLD}}},
      \bQ_{t,\delta}\right)
 	  \le \big(\phi(t) + {\mu_2^{1/2}(\pi)}\big)
 	  \exp\bigg\{M\delta  +  \int_0^\delta \alpha(t+u)\,du \bigg\}
 	  \int_0^\delta \big|\alpha(t+s)-\alpha(t)\big|\,ds.
   \end{equation}
\end{lemma}   

When $\delta$ tends to $0$, according to \Cref{lem-gr-1st}, the first term of the 
right-hand side of \eqref{eq-minus-phi} vanishes. Thus, after passing to the limit, we 
are left with the  following  Gronwall-type	inequality:
\begin{equation}\label{eq-gron-phi}
	\phi'(t)\leq  -(m+\alpha(t)) \phi(t)  -
	\frac{11\alpha'(t)}{\sqrt{m+\alpha(t)}}\cdot 
	\mu_2(\pi_{\alpha(t)}).
\end{equation}
Here we tacitly used the fact that $\mu_2(\pi_{\alpha(t+\delta)}) \rightarrow
\mu_2 ( \pi_{\alpha(t)})$, whenever $\delta \rightarrow 0$, which is due to the 
continuity of $\alpha(t)$ and \Cref{lem-mon-mu}. Recalling that the function 
$\beta(t)$ is given by $\beta(t) = \int_0^t\big(m+\alpha(s)\big)\,ds$, one can 
rewrite  
 \eqref{eq-gron-phi}  as
\begin{equation}
	\big(\phi(t) e^{\beta(t)} \big)'
	\leq -  \frac{11\alpha'(t) e^{\beta(t)}}{\sqrt{m+\alpha(t)}} \mu_2(\pi_{\alpha(t)})
	\leq -  \frac{11\alpha'(t) e^{\beta(t)}}{\sqrt{m+\alpha(t)}} \mu_2(\pi).
\end{equation}
Therefore we infer the following bound on $\phi(t)$:
\begin{align}\label{eq-before-ipp}
 	\phi(t)&\leq  \phi(0) e^{-\beta(t)} -
	11\mu_2(\pi)\int_0^t \frac{\alpha'(s)}{\sqrt{m+\alpha(s)}} 
	e^{\beta(s)-\beta(t)} ds.
\end{align}
Combining this bound with \eqref{eq-first-tr-ineq} and  \eqref{eq-last-term}, we obtain the inequality
\begin{align}
	W_2(\nu_{t}^{\textsf{PLD}},\pi)
	& \leq  W_2(\nu_0,\pi_{\alpha(0)}) e^{-\beta(t)}  - 11\mu_2(\pi)\int_0^t \frac{\alpha'(s)e^{\beta(s)-\beta(t)}}{\sqrt{m+\alpha(s)}} ds+
	\frac{11\alpha(t)\mu_2(\pi)}{\sqrt{m+\alpha(t)}}.
\end{align}
 \Cref{lem-mon-mu} yields
$W_2(\nu_0,\pi_{\alpha(0)}) = \sqrt{\mu_2(\pi_{\alpha(0)})} \leq \sqrt{\mu_2(\pi)}$. This completes the proof of \Cref{th-main}, since the derivative of $\alpha$ is  negative.


\bibliography{Literature1}

\begin{thebibliography}{}

\bibitem[\protect\astroncite{Alonso-Guti\'{e}rrez and
  Bastero}{2015}]{KLS_conjecture}
Alonso-Guti\'{e}rrez, D. and Bastero, J. (2015).
\newblock {\em Approaching the {K}annan-{L}ov\'{a}sz-{S}imonovits and variance
  conjectures}, volume 2131 of {\em Lecture Notes in Mathematics}.
\newblock Springer, Cham.

\bibitem[\protect\astroncite{Ambrosio et~al.}{2008}]{Ambrosio}
Ambrosio, L., Gigli, N., and Savar\'{e}, G. (2008).
\newblock {\em Gradient flows in metric spaces and in the space of probability
  measures}.
\newblock Lectures in Mathematics ETH Z\"{u}rich. Birkh\"{a}user Verlag, Basel,
  second edition.

\bibitem[\protect\astroncite{Arbel et~al.}{2019}]{Arbel19}
Arbel, M., Korba, A., Salim, A., and Gretton, A. (2019).
\newblock Maximum mean discrepancy gradient flow.
\newblock In {\em Advances in Neural Information Processing Systems 32}, pages
  6484--6494. Curran Associates, Inc.

\bibitem[\protect\astroncite{Bernton}{2018}]{Bernton18}
Bernton, E. (2018).
\newblock Langevin monte carlo and {JKO} splitting.
\newblock In {\em Conference On Learning Theory, {COLT} 2018, Stockholm,
  Sweden, 6-9 July 2018}, volume~75 of {\em Proceedings of Machine Learning
  Research}, pages 1777--1798. {PMLR}.

\bibitem[\protect\astroncite{Bhattacharya}{1978}]{bhattacharya1978}
Bhattacharya, R.~N. (1978).
\newblock Criteria for recurrence and existence of invariant measures for
  multidimensional diffusions.
\newblock {\em Ann. Probab.}, 6(4):541--553.

\bibitem[\protect\astroncite{Bobkov}{1999}]{Bobkov99}
Bobkov, S.~G. (1999).
\newblock Isoperimetric and analytic inequalities for log-concave probability
  measures.
\newblock {\em Ann. Probab.}, 27(4):1903--1921.

\bibitem[\protect\astroncite{Bolley et~al.}{2012}]{Bolley}
Bolley, F., Gentil, I., and Guillin, A. (2012).
\newblock Convergence to equilibrium in {W}asserstein distance for
  {F}okker-{P}lanck equations.
\newblock {\em J. Funct. Anal.}, 263(8):2430--2457.

\bibitem[\protect\astroncite{Bubeck et~al.}{2018}]{Bubeck18}
Bubeck, S., Eldan, R., and Lehec, J. (2018).
\newblock Sampling from a log-concave distribution with projected langevin
  monte carlo.
\newblock {\em Discrete {\&} Computational Geometry}, 59(4):757--783.

\bibitem[\protect\astroncite{Cattiaux and Guillin}{2009}]{Cattiaux1}
Cattiaux, P. and Guillin, A. (2009).
\newblock Trends to equilibrium in total variation distance.
\newblock {\em Ann. Inst. Henri Poincar\'{e} Probab. Stat.}, 45(1):117--145.

\bibitem[\protect\astroncite{Cattiaux and Guillin}{2018}]{cattiaux2018poincar}
Cattiaux, P. and Guillin, A. (2018).
\newblock On the {P}oincar{\'e} constant of log-concave measures.
\newblock {\em arXiv preprint arXiv:1810.08369}.

\bibitem[\protect\astroncite{Cattiaux et~al.}{2010}]{Cattiaux2}
Cattiaux, P., Guillin, A., and Roberto, C. (2010).
\newblock Poincar\'{e} inequality and the {$L^p$} convergence of semi-groups.
\newblock {\em Electron. Commun. Probab.}, 15:270--280.

\bibitem[\protect\astroncite{Chatterji et~al.}{2019}]{chatterji2019langevin}
Chatterji, N.~S., Diakonikolas, J., Jordan, M.~I., and Bartlett, P.~L. (2019).
\newblock Langevin monte carlo without smoothness.
\newblock {\em arXiv preprint arXiv:1905.13285}.

\bibitem[\protect\astroncite{Cheng and Bartlett}{2018}]{Cheng1}
Cheng, X. and Bartlett, P. (2018).
\newblock Convergence of {L}angevin {MCMC} in {KL}-divergence.
\newblock In {\em Proceedings of ALT2018}.

\bibitem[\protect\astroncite{Cheng et~al.}{2018}]{Cheng3}
Cheng, X., Chatterji, N.~S., Abbasi{-}Yadkori, Y., Bartlett, P.~L., and Jordan,
  M.~I. (2018).
\newblock Sharp convergence rates for langevin dynamics in the nonconvex
  setting.
\newblock {\em CoRR}, abs/1805.01648.

\bibitem[\protect\astroncite{{Cheng} et~al.}{2017}]{Cheng2}
{Cheng}, X., {Chatterji}, N.~S., {Bartlett}, P.~L., and {Jordan}, M.~I. (2017).
\newblock {Underdamped Langevin MCMC: A non-asymptotic analysis}.
\newblock {\em ArXiv e-prints}.

\bibitem[\protect\astroncite{Chewi et~al.}{2020}]{chewi2020exponential}
Chewi, S., Gouic, T.~L., Lu, C., Maunu, T., Rigollet, P., and Stromme, A.
  (2020).
\newblock Exponential ergodicity of mirror-{L}angevin diffusions.

\bibitem[\protect\astroncite{Dalalyan}{2017}]{Dalalyan14}
Dalalyan, A.~S. (2017).
\newblock Theoretical guarantees for approximate sampling from a smooth and
  log-concave density.
\newblock {\em J. R. Stat. Soc. B}, 79:651 -- 676.

\bibitem[\protect\astroncite{Dalalyan and
  Riou-Durand}{2018}]{dalalyan_riou_2018}
Dalalyan, A.~S. and Riou-Durand, L. (2018).
\newblock On sampling from a log-concave density using kinetic langevin
  diffusions.
\newblock {\em arXiv preprint arXiv:1807.09382}.

\bibitem[\protect\astroncite{Dalalyan et~al.}{2019}]{dalalyan2019bounding}
Dalalyan, A.~S., Riou-Durand, L., and Karagulyan, A. (2019).
\newblock Bounding the error of discretized langevin algorithms for
  non-strongly log-concave targets.

\bibitem[\protect\astroncite{Durmus et~al.}{2018a}]{durmus2018analysis}
Durmus, A., Majewski, S., and Miasojedow, B. (2018a).
\newblock Analysis of langevin monte carlo via convex optimization.
\newblock {\em arXiv preprint arXiv:1802.09188}.

\bibitem[\protect\astroncite{Durmus and Moulines}{2017}]{durmus2017}
Durmus, A. and Moulines, E. (2017).
\newblock Nonasymptotic convergence analysis for the unadjusted {L}angevin
  algorithm.
\newblock {\em Ann. Appl. Probab.}, 27(3):1551--1587.

\bibitem[\protect\astroncite{Durmus and Moulines}{2019}]{Durmus2}
Durmus, A. and Moulines, E. (2019).
\newblock High-dimensional bayesian inference via the unadjusted langevin
  algorithm.
\newblock {\em Bernoulli}, 25(4A):2854--2882.

\bibitem[\protect\astroncite{Durmus et~al.}{2018b}]{Durmus3}
Durmus, A., Moulines, {\'E}., and Pereyra, M. (2018b).
\newblock {Efficient Bayesian Computation by Proximal Markov Chain Monte Carlo:
  When Langevin Meets Moreau}.
\newblock {\em {SIAM Journal on Imaging Sciences}}, 11(1).

\bibitem[\protect\astroncite{{Dwivedi} et~al.}{2018}]{Dwivedi}
{Dwivedi}, R., {Chen}, Y., {Wainwright}, M.~J., and {Yu}, B. (2018).
\newblock {Log-concave sampling: Metropolis-Hastings algorithms are fast}.
\newblock {\em arXiv e-prints}.

\bibitem[\protect\astroncite{Eberle et~al.}{2019}]{eberle2019}
Eberle, A., Guillin, A., and Zimmer, R. (2019).
\newblock Couplings and quantitative contraction rates for langevin dynamics.
\newblock {\em Ann. Probab.}, 47(4):1982--2010.

\bibitem[\protect\astroncite{Erdogdu et~al.}{2018}]{Erdogdu18}
Erdogdu, M.~A., Mackey, L., and Shamir, O. (2018).
\newblock Global non-convex optimization with discretized diffusions.
\newblock In {\em Advances in Neural Information Processing Systems 31}, pages
  9671--9680.

\bibitem[\protect\astroncite{Franca et~al.}{2018}]{franca18a}
Franca, G., Robinson, D., and Vidal, R. (2018).
\newblock {ADMM} and accelerated {ADMM} as continuous dynamical systems.
\newblock In {\em Proceedings of the 35th International Conference on Machine
  Learning}, volume~80 of {\em Proceedings of Machine Learning Research}, pages
  1559--1567. PMLR.

\bibitem[\protect\astroncite{Gozlan and
  L\'{e}onard}{2010}]{gozlan2010transport}
Gozlan, N. and L\'{e}onard, C. (2010).
\newblock Transport inequalities. {A} survey.
\newblock {\em Markov Process. Related Fields}, 16(4):635--736.

\bibitem[\protect\astroncite{Hsieh et~al.}{2018}]{Hsieh:266354}
Hsieh, Y.-P., Kavis, A., Rolland, P., and Cevher, V. (2018).
\newblock Mirrored langevin dynamics.
\newblock {\em Advances In Neural Information Processing Systems 31 (Nips
  2018)}, 31.

\bibitem[\protect\astroncite{Kannan et~al.}{1995}]{Kannan95}
Kannan, R., Lov\'{a}sz, L., and Simonovits, M. (1995).
\newblock Isoperimetric problems for convex bodies and a localization lemma.
\newblock {\em Discrete Comput. Geom.}, 13(3-4):541--559.

\bibitem[\protect\astroncite{Krichene et~al.}{2015}]{Krichene15}
Krichene, W., Bayen, A., and Bartlett, P.~L. (2015).
\newblock Accelerated mirror descent in continuous and discrete time.
\newblock In Cortes, C., Lawrence, N.~D., Lee, D.~D., Sugiyama, M., and
  Garnett, R., editors, {\em Advances in Neural Information Processing Systems
  28}, pages 2845--2853. Curran Associates, Inc.

\bibitem[\protect\astroncite{Liang et~al.}{2019}]{liang2019exponential}
Liang, M., Majka, M.~B., and Wang, J. (2019).
\newblock Exponential ergodicity for sdes and mckean-vlasov processes with
  lévy noise.

\bibitem[\protect\astroncite{Ma et~al.}{2019}]{YiAnMa}
Ma, Y., Chatterji, N.~S., Cheng, X., Flammarion, N., Bartlett, P.~L., and
  Jordan, M.~I. (2019).
\newblock Is there an analog of nesterov acceleration for mcmc?
\newblock {\em CoRR}, abs/1902.00996.

\bibitem[\protect\astroncite{Majka et~al.}{2018}]{majka2018nonasymptotic}
Majka, M.~B., Mijatović, A., and Szpruch, L. (2018).
\newblock Non-asymptotic bounds for sampling algorithms without log-concavity.

\bibitem[\protect\astroncite{Mangoubi and Vishnoi}{2019}]{MangoubiV19}
Mangoubi, O. and Vishnoi, N.~K. (2019).
\newblock Nonconvex sampling with the metropolis-adjusted langevin algorithm.
\newblock In Beygelzimer, A. and Hsu, D., editors, {\em Conference on Learning
  Theory, {COLT} 2019, 25-28 June 2019, Phoenix, AZ, {USA}}, volume~99 of {\em
  Proceedings of Machine Learning Research}, pages 2259--2293. {PMLR}.

\bibitem[\protect\astroncite{Mou et~al.}{2019}]{Mou}
Mou, W., Flammarion, N., Wainwright, M.~J., and Bartlett, P.~L. (2019).
\newblock An efficient sampling algorithm for non-smooth composite potentials.
\newblock {\em CoRR}, abs/1910.00551.

\bibitem[\protect\astroncite{Nesterov}{2004}]{Nest}
Nesterov, Y. (2004).
\newblock {\em Introductory lectures on convex optimization}, volume~87 of {\em
  Applied Optimization}.
\newblock Kluwer Academic Publishers, Boston, MA.

\bibitem[\protect\astroncite{Pillaud-Vivien
  et~al.}{2019}]{pillaudvivien2019statistical}
Pillaud-Vivien, L., Bach, F., Lelièvre, T., Rudi, A., and Stoltz, G. (2019).
\newblock Statistical estimation of the poincar{é} constant and application to
  sampling multimodal distributions.

\bibitem[\protect\astroncite{Salim et~al.}{2019}]{Salim}
Salim, A., Koralev, D., and Richtarik, P. (2019).
\newblock Stochastic proximal langevin algorithm: Potential splitting and
  nonasymptotic rates.
\newblock In {\em Advances in Neural Information Processing Systems 32}, pages
  6653--6664.

\bibitem[\protect\astroncite{Santambrogio}{2017}]{Santa17}
Santambrogio, F. (2017).
\newblock \{{E}uclidean, metric, and {W}asserstein\} gradient flows: an
  overview.
\newblock {\em Bull. Math. Sci.}, 7(1):87--154.

\bibitem[\protect\astroncite{Scieur et~al.}{2017}]{Scieur}
Scieur, D., Roulet, V., Bach, F., and d'Aspremont, A. (2017).
\newblock Integration methods and optimization algorithms.
\newblock In {\em Advances in Neural Information Processing Systems 30}, pages
  1109--1118. Curran Associates, Inc.

\bibitem[\protect\astroncite{Shen and Lee}{2019}]{shen2019randomized}
Shen, R. and Lee, Y.~T. (2019).
\newblock The randomized midpoint method for log-concave sampling.
\newblock In {\em Advances in Neural Information Processing Systems}, pages
  2098--2109.

\bibitem[\protect\astroncite{Simsekli et~al.}{2020}]{umut20}
Simsekli, U., Zhu, L., Teh, Y.~W., and G{\"{u}}rb{\"{u}}zbalaban, M. (2020).
\newblock Fractional underdamped langevin dynamics: Retargeting {SGD} with
  momentum under heavy-tailed gradient noise.
\newblock {\em CoRR}, abs/2002.05685.

\bibitem[\protect\astroncite{Su et~al.}{2016}]{SuBoydCandes}
Su, W., Boyd, S., and Cand{{\`e}}s, E.~J. (2016).
\newblock A differential equation for modeling nesterov's accelerated gradient
  method: Theory and insights.
\newblock {\em Journal of Machine Learning Research}, 17(153):1--43.

\bibitem[\protect\astroncite{Vempala and Wibisono}{2019}]{Vempala_Wibisono}
Vempala, S. and Wibisono, A. (2019).
\newblock Rapid convergence of the unadjusted langevin algorithm: Isoperimetry
  suffices.
\newblock In {\em Advances in Neural Information Processing Systems 32}, pages
  8094--8106.

\bibitem[\protect\astroncite{Villani}{2008}]{villani2008optimal}
Villani, C. (2008).
\newblock {\em Optimal transport: old and new}, volume 338.
\newblock Springer Science \& Business Media.

\bibitem[\protect\astroncite{Wibisono}{2018}]{wibisono18a}
Wibisono, A. (2018).
\newblock Sampling as optimization in the space of measures: The langevin
  dynamics as a composite optimization problem.
\newblock In {\em Proceedings of the 31st Conference On Learning Theory},
  volume~75 of {\em Proceedings of Machine Learning Research}, pages
  2093--3027. PMLR.

\bibitem[\protect\astroncite{Wibisono et~al.}{2016}]{Wibisono}
Wibisono, A., Wilson, A.~C., and Jordan, M.~I. (2016).
\newblock A variational perspective on accelerated methods in optimization.
\newblock {\em Proceedings of the National Academy of Sciences},
  113(47):E7351--E7358.

\bibitem[\protect\astroncite{Wilson et~al.}{2016}]{wilson2016lyapunov}
Wilson, A.~C., Recht, B., and Jordan, M.~I. (2016).
\newblock A lyapunov analysis of momentum methods in optimization.

\bibitem[\protect\astroncite{Zhang et~al.}{2018}]{Jingzhao18}
Zhang, J., Mokhtari, A., Sra, S., and Jadbabaie, A. (2018).
\newblock Direct runge-kutta discretization achieves acceleration.
\newblock In Bengio, S., Wallach, H., Larochelle, H., Grauman, K.,
  Cesa-Bianchi, N., and Garnett, R., editors, {\em Advances in Neural
  Information Processing Systems 31}, pages 3900--3909.

\end{thebibliography}

\newpage
\appendix	
\renewcommand\contentsname{}
\part{Appendix} 
\parttoc 


\section{Proofs of the lemmas used in \Cref{th-main}}\label{sec-appendix}

\subsection{Proof of \Cref{lem-tr-ineq}}

We denote by $D_{\rm KL}(\pi_{\gamma}||\pi_{{\tilde{\gamma}}})$
the Kullback-Leibler divergence between the distributions
$\pi_{\gamma}$ and $\pi_{\tilde{\gamma}}$. Since 
$\pi_{\tilde{\gamma}}$ is $m+\tilde\gamma$-strongly log-concave,
the transportation cost inequality \citep[Corollary 7.2]
{gozlan2010transport} yields
\begin{equation}\label{gozlan-tr-ineq}
    W_2^2(\pi_{{\tilde{\gamma}}},\pi_{\gamma})
        \le \frac{2}{m+{\tilde{\gamma}}} D_{\rm KL}(\pi_{\gamma}||\pi_{{\tilde{\gamma}}}).
\end{equation}
Let us denote by $c_\gamma$ the logarithm of the normalizing
constant for $\pi_{\gamma}$ so that $\pi_{\gamma}(\btheta) =
\exp(-f(\btheta) - (\nicefrac12)\gamma \| \btheta \|_2^2 +
c_\gamma)$. Similarly, we denote by $c_{{\tilde{\gamma}}}$ the
logarithm of the normalizing constant of $\pi_{\tilde\gamma}$. 
This readily yields
\begin{align}
    D_{\rm KL}(\pi_{\gamma}||\pi_{{\tilde{\gamma}}})
    &= \int_{\mathbb{R}^p} \pi_\gamma(\btheta) 
    \log\left( \frac{\pi_\gamma(\btheta)}{\pi_{\tilde{\gamma}}(\btheta)}\right)d\btheta \\
    &= \int_{\mathbb{R}^p} \pi_\gamma(\btheta) \left(\nicefrac12({\tilde{\gamma}} - \gamma)
    \|\btheta\|_2^2 +  c_\gamma - c_{\tilde{\gamma}}\right)d\btheta\\
    &= \nicefrac12({\tilde{\gamma}} - \gamma)\mu _2(\pi_\gamma) + c_\gamma - c_{\tilde{\gamma}}.
\end{align}
Using the inequality $e^{-x}\le 1-x + (\nicefrac12)x^2$ for all $x>0$
implies the following upper bound on $c_\gamma - c_{\tilde\gamma}$:
\begin{align}
	c_\gamma - c_{\tilde{\gamma}}  &= \log\left( \int_{\RR^p} \pi_\gamma(\btheta) 
	\exp\left(\nicefrac12(\gamma - {\tilde{\gamma}})\|\btheta\|_2^2 \right)\right)	\\
	&\leq  \log\left( 1 + \nicefrac{1}{2}(\gamma-{\tilde{\gamma}}) \mu_2(\pi_\gamma)  +   
    \nicefrac{1}{8}(\gamma - {\tilde{\gamma}})^2\mu_4(\pi_\gamma) \right).
\end{align}	    	
Since $\log(1+x)\le x$ for $x>-1$ we get
\begin{equation}
 D_{\rm KL}(\pi_{\gamma}||\pi_{{\tilde{\gamma}}}) \le 
 \nicefrac{1}{8} (\gamma - \tilde{\gamma})^2\mu_4(\pi_\gamma).
\end{equation}
Since $m+\gamma\geq 0$, the distribution $\pi_\gamma$ is
log-concave. Thus, in view of \citep[Remark 3]{dalalyan2019bounding},
we have the inequality $\mu_{4}(\pi_\gamma) \leq 442 
\mu_{2}^2(\pi_\gamma)$. Finally, combining these bounds with
\eqref{gozlan-tr-ineq}, we get
\begin{equation}
    W_2(\pi_{\tilde{\gamma}},\pi_\gamma)
        \le \sqrt{\frac{2}{m+\tilde\gamma}} \times \frac{(\tilde{\gamma} - \gamma) \mu_4^{1/2}(\pi_\gamma)}
        {\sqrt{8}} 
        \le  \frac{11\mu_2(\pi_\gamma)}{\sqrt{m+{\tilde{\gamma}}}}
        ({\tilde{\gamma}} - \gamma).
\end{equation}
This completes the proof of the lemma.

\subsection{Proof of \Cref{lem-mon-mu}}

	For $k\in \NN \cup \{0\}$, define  
	\begin{equation}
		h_k(\gamma) = \int_{\RR^p} \|\btheta\|_2^k \exp \left(-
		f(\btheta) - {\gamma\|\btheta\|_2^2}/2\right)d\btheta.
	\end{equation}
	If $\pi \in \mathcal{P}_{k} (\RR^p)$ then the function $h_k$ is
	continuous on $[0;+\infty)$. Indeed, if the sequence
	$\{\gamma_n\}_{n}$ converges $\gamma_0$, when $n\rightarrow +\infty$, 
	then the function  $\|\btheta\|_2^k \exp \left(-f(\btheta) - 	(\nicefrac{1}{2}){\gamma_n \|\btheta\|_2^2}\right)$ is 
	upper-bounded by $ \|\btheta\|_2^k \exp \left(- f(\btheta)\right)$.  Thus in view of the dominated convergence theorem, we can interchange the limit and the integral. Since, by definition, 
	\begin{equation}
		\mu_k(\pi_{\gamma}) = \frac{h_k(\gamma)}{h_0(\gamma)},
	\end{equation}
	we get the continuity of $\mu_2(\pi_{\gamma})$ and $\mu_4(\pi_{\gamma})$. 
	Let us now prove that $h_k(t)$ is continuously differentiable,
	 when $\pi \in \mathcal{P}_{k+2} (\RR^p) $.  The integrand function in the definition of $h_k$
	 is a continuously differentiable function with respect to $t$. In addition, its  derivative is   
	 continuous and is as well integrable on $\RR^p$,  as we supposed that $\pi$ 
	 has the $(k+2)$-th 
	 moment.  Therefore, Leibniz integral rule yields the following 
	 \begin{equation}
	 	h_k'(\gamma) = -\frac{1}{2}\int_{\RR^p} \|\btheta\|_2^{k+2} \exp \left(-
		f(\btheta) - {\gamma \|\btheta\|_2^2}/{2}\right)d\btheta =  -\frac{1}{2}h_{k+2}(t).
	 \end{equation}
    The latter yields the smoothness of $h_k$. 
	Finally, in order to  prove the monotony of $\mu_2(\pi_{\gamma})$, we will
	simply calculate its derivative 
	\begin{align}
		\left(\mu_2(\pi_{\gamma})\right)' &= - \frac{1}{2h_{0}(\gamma)} 
		h_4(\gamma)	- \frac{h_{0}'(\gamma)}{h_{0} (\gamma)^2}{h_2}(\gamma)\\
		&= -\frac{1}{2}\mu_4(\pi_{\gamma}) +  
		\frac{h_{2}^2(\gamma)}{2h_{0}(\gamma)^2} \\
		&=  \frac{1}{2}\left( \mu_2^2(\pi_{\gamma}) - \mu_4(\pi_{\gamma})\right).
	\end{align}
	Since the latter is always negative, this completes the proof of the lemma.

\subsection{Proof of \Cref{lem-gr-1st}}
	From the definition of Wasserstein distance, we have
    \begin{equation}
        W_2\left( \nu_{t+\delta}^{\textupsf{PLD}},\bQ_{t,\delta} \right) \leq \|\tilde{\bL}_{t+\delta} - \bL_{t+\delta}^{\textupsf{PLD}}\|_{\mathbb L_2}.
    \end{equation}
 	In view of the definition of the process $\tilde{\bL}$, we 
    can write
    \begin{equation}
        \tilde{\bL}_{t+\delta} - \bL_{t+\delta}^{\textupsf{PLD}} 
        = \int_{t}^{t+\delta} \left(\nabla f({\bL}_{s}^{\textupsf{PLD}}) - 
        \nabla f(\tilde{\bL}_{s}^{\textupsf{PLD}}) + 
        \alpha(s){\bL}_{s}^{\textupsf{PLD}} - \alpha(t)\tilde{\bL}_{s}\right)ds.
    \end{equation}
    Therefore we have
    \begin{align}
        \|\tilde{\bL}_{t+\delta} - \bL_{t+\delta}^{\textupsf{PLD}}
        \|_{\mathbb L_2} \leq  &
        \bigg\| \underbrace{\int_{t}^{t+\delta} \left(\nabla f({\bL}_{s}^{\textupsf{PLD}}) - 
        \nabla f(\tilde{\bL}_{s})\right)ds}_{:=T_1}\bigg\|_{\mathbb L_2} + \bigg\| \underbrace{\int_{t}^{t+\delta} \left(\alpha(s)
        {\bL}_{s}^{\textupsf{PLD}} - \alpha(t)\tilde{\bL}_{s}\right)ds}_{:=T_2}\bigg\|_{\mathbb L_2}. 
        \label{eq:11}
    \end{align}    Now let us analyze these two terms separately. We start with  $T_1$:
    \begin{align}
         \|T_1\|_{\mathbb L_2} &= \left\| \int_{t}^{t+\delta} \left(\nabla f({\bL}_{s}^{\textupsf{PLD}}) - 
        \nabla f(\tilde{\bL}_{s})\right)ds\right\|_{\mathbb L_2}\\ &\leq
         \int_{t}^{t+\delta} \left\|\nabla f({\bL}_{s}^{\textupsf{PLD}}) - 
        \nabla f(\tilde{\bL}_{s})\right\|_{\mathbb L_2}ds \\
        & \leq M\int_{t}^{t+\delta}  \|{\bL}_{s}^{\textupsf{PLD}} - 
        \tilde{\bL}_{s}\|_{\mathbb L_2}ds.
    \end{align}
    These are due to the Minkowskii inequality and the Lipschitz 
    continuity of the gradient.
    In order to bound the second term $T_2$, we will add and subtract 
    the term 
    $\alpha(t+s)\tilde{\bL}_{t+s}$.  Similar to the case above,
     we get the 
    following upper bound:
    \begin{align}
    	\|T_2\|_{\mathbb L_2} &\le  \int_{t}^{t+\delta} \alpha(s) 
         \,\big\| {\bL}_{s}^{\textupsf{PLD}} - \tilde{\bL}_{s}\big\|_{\mathbb L_2} ds 
         +  
         \int_{t}^{t+\delta} \big|\alpha(s) - \alpha(t)\big|\,\big\|\tilde{\bL}_{s}\big\|_{\mathbb L_2}
          ds\\
          &=\int_{0}^{\delta} \alpha(t+s) 
         \,\big\| {\bL}_{t+s}^{\textupsf{PLD}} - \tilde{\bL}_{t+s}\big\|_{\mathbb L_2} ds 
         +  
         \int_{0}^{\delta} \big|\alpha(t+s) 
         - \alpha(t)\big|\,\big\|\tilde{\bL}_{t+s} 
         \big\|_{\mathbb L_2} ds.
    \end{align}
    Recall that 
    $\tilde{\bL}_{t+s}$ is the solution of Langevin SDE with an 
    $(m+\alpha(t))$-strongly convex  potential function, and $Q_{t,s}$
    is its distribution on $\mathbb R^p$. Thus, the triangle inequality
    yields 
    \begin{align}
    	\big\|\tilde{\bL}_{t+s}\big\|_{\mathbb L_2}  
    	    &= W_2(\bQ_{t,s},\delta_0) 
    	    \le W_2(\bQ_{t,s},\pi_{\alpha(t)}) + W_2(\pi_{\alpha(t)},\delta_0)\\
    		&\leq W_2(\nu_t^{\textupsf{PLD}},\pi_{\alpha(t)})\exp(-ms-\alpha(t)s) + \sqrt{\mu _2(\pi_{\alpha(t)})}\\ 
    		&\leq   W_2(\nu^{\textupsf{PLD}}_t,\pi_{\alpha(t)}) + \sqrt{\mu_2(\pi_{\alpha(t)})}:=V_t.
    \end{align}
    Summing up, we have
    \begin{align}
         \big\|  {\bL}_{t+\delta}^{\textupsf{PLD}} - \tilde{\bL}_{t+\delta}
         \big\|_{\mathbb L_2} &\leq
          \int_{0}^\delta \big(M +  \alpha(t+s)\big)
         \|  {\bL}_{t+s}^{\textupsf{PLD}} - \tilde{\bL}_{t+s}\|_{\mathbb L_2}
         ds + \tilde\alpha_t(\delta)\,V_t,
	\end{align}
	where $\tilde{\alpha_t}(\delta)$ is an auxiliary function defined as
	\begin{align}
        \tilde \alpha_t(\delta) :=  \int_0^\delta
        |\alpha(t+s)-\alpha(t)|\,ds.
    \end{align}
Now let us define  $\Phi(s) =  \|  {\bL}_{t+s}^{\textupsf{PLD}} - \tilde{\bL}_{t+s}\|_{\mathbb L_2}$. The last inequality can be rewritten as
	\begin{equation}
		\Phi(\delta) \leq \int_0^\delta \big(M  +  \alpha(t+s)\big)\Phi(s)\,ds
		+\tilde\alpha_t(\delta)\, V_t.
	\end{equation}
The (integral form of the) Gronwall inequality, \cref{lem:Gron_int}, implies that
\begin{align}
    \Phi(\delta)  &\le V_t
    \int_0^\delta 
    \tilde\alpha_t(s)\big(M  +  \alpha(t+s)\big) e^{\int_s^\delta (M  +  \alpha(t+u))\,du}\,ds + \tilde\alpha_t(\delta)V_t \\
    &= V_t \int_0^\delta 
    \tilde\alpha_t'(s)\,e^{\int_s^\delta (M  +  \alpha(t+u))\,du}\,ds\\
    &\le V_t\,\tilde\alpha_t(\delta)\,
    \exp\bigg\{M\delta  +  \int_0^\delta \alpha(t+u)\,du \bigg\}. 
\end{align}
This completes the proof. 

\subsection{Different forms of the Gronwall inequality}

In this section, we provide two forms of the Gronwall inequality
that are used in the present work. For the sake of the self-containedness, the proofs of these inequalities are also 
provided.

\begin{lemma}[Differential form]\label{lem:Gron_diff}
Let $A:[a,b]\to \mathbb R$ and $B:[a,b]\to \mathbb R$ be two 
functions. If the function $\Phi:[a,b]\to \mathbb R$ satisfies
the recursive differential inequality 
\begin{align}
    \Phi'(x) \le A(x)\Phi(x) + B(x),\qquad 
    \forall x\in[a,b],
\end{align}
then it also satisfies the inequality
\begin{align}
    \Phi(x) \le \Phi(a)\exp\bigg\{\int_a^x A(z)\,dz\bigg\} + 
    \int_a^x B(s) \exp\bigg\{\int_s^x A(z)\,dz\bigg\}\,ds,\qquad 
    \forall x\in[a,b].
\end{align}
\end{lemma}

\begin{proof}
To be inserted. 
\end{proof}

\begin{lemma}[Integral form]\label{lem:Gron_int}
Let $A:[a,b]\to [0,+\infty)$ and $B:[a,b]\to \mathbb R$ be two 
functions. If the function $\Phi:[a,b]\to \mathbb R$ satisfies
the recursive integral inequality 
\begin{align}
    \Phi(x) \le \int_a^x A(s)\Phi(s)\,ds + B(x),\qquad 
    \forall x\in[a,b],
\end{align}
then it also satisfies the inequality
\begin{align}\label{eq:Phi2}
    \Phi(x) \le \int_a^x A(s)B(s) \exp\bigg\{\int_s^x A(z)\,dz\bigg\}\,ds + B(x),\qquad 
    \forall x\in[a,b].
\end{align}
\end{lemma}

\begin{proof}
To be inserted. We set 
\begin{align}
    \Psi(x) = \exp\bigg\{-\int_a^x A(z)\,dz\bigg\} 
    \int_a^x A(s)\Phi(s)\,ds. 
\end{align}
We have
\begin{align}
    \Psi'(x) &= -A(x)\Psi(x) +\exp\bigg\{-\int_a^x A(z)\,dz\bigg\} 
     A(x)\Phi(x)\\
     &\le -A(x)\Psi(x) +\exp\bigg\{-\int_a^x A(z)\,dz\bigg\} 
     A(x)\Big(\int_a^x A(s)\Phi(s)\,ds + B(x)\Big)\\
     & = -A(x)\Psi(x) + 
     A(x)\Psi(x) + A(x)B(x)\exp\bigg\{-\int_a^x A(z)\,dz\bigg\}.
\end{align}
Therefore,
\begin{align}
    \Psi(x) &\le \Psi(a) + \int_a^x A(s)B(s)\exp\bigg\{-\int_a^s A(z)\,dz\bigg\}\,ds.
\end{align}
Replacing $\Psi$ by its expression and using the fact that $\Psi(a) = 0$, we get 
\begin{align}
    \exp\bigg\{-\int_a^x A(z)\,dz\bigg\} 
    \int_a^x A(s)\Phi(s)\,ds &\le \int_a^x A(s)B(s)\exp\bigg\{-\int_a^s A(z)\,dz\bigg\}\,ds.
\end{align}
This implies that
\begin{align}
    \int_a^x A(s)\Phi(s)\,ds &\le \int_a^x A(s)B(s)\exp\bigg\{\int_s^x A(z)\,dz\bigg\}\,ds.
\end{align}
Combining this inequality with \eqref{eq:Phi2}, we get the claim 
of the lemma. 
\end{proof}

\section{Proof of \protect\Cref{prop-optimal}}

For the penalty factor $\alpha(t) = 1/({A+2t})$, we get $\beta(t) =
\int_0^t \alpha(s)\,ds = (\nicefrac12) \log\big(1+(2/A)t\big)$. This
implies that
\begin{align}
    \sqrt{\mu_2(\pi)}\, e^{-\beta(t)}  + 11\mu_2(\pi)\sqrt{\alpha(t)} = 
    \frac{\sqrt{A\mu_2(\pi)} + 11 \mu_2(\pi)}{\sqrt{A+2t}}. 
\end{align}
Finally, the middle term in the right hand side of \eqref{main:ineq}
takes the form
\begin{align}
    11\mu_2(\pi)\int_0^t \frac{|\alpha'(s)|}{\sqrt{\alpha(s)}} e^{\beta(s)- \beta(t)} ds 
    & = \frac{11\mu_2(\pi)}{\sqrt{A + 2t}}\int_0^t \frac{2}{A + 2s} ds\\
    & = \frac{11\mu_2(\pi)}{\sqrt{A + 2t}} \log{\left(1 + (2/A)t\right)}.
\end{align}
Combining these relations, we get the claim of the proposition.

\section{(Weakly) convex potentials: what is known and what we 
can hope for} 

Many recent papers investigated the case of strongly convex potential;
this case is now rather well understood. Let us briefly summarize here
some facts and conjectures that can shed some light on the broader
case of weakly convex potential. This might help to understand what 
can be expected to be proved in the framework s<tudied in this work. 

The ergodicity properties of the Langevin process are closely related
to such notions of functional analysis as the spectral gap, the Poincar\'e and the log-Sobolev inequalities. Thus, the generator of a Markov semi-group associated with an $m$-strongly convex potential has a spectral-gap $\mathcal C_{\textup{SG}}$ at least equal to $m$. This property was exploited by \cite{Dalalyan14} to derive guarantees on the LMC algorithm. It is known that the spectral gap exists if and only if the invariant density satisfies the Poincar\'e inequality. Furthermore, the spectral gap is equal to the inverse of
the Poincar\'e constant $\mathcal C_{\textup{P}}$. Furthermore, distributions associated to $m$-strongly convex potentials satisfy 
the log-Sobolev inequality with the constant $\mathcal C_{\textup{LS}}\le 1/m$. This property was used by \cite{Durmus2} 
to extend the guarantees to the Wasserstein-2 distance. 

Note that the log-Sobolev inequality is stronger than the Poincar\'e inequality and  $\mathcal C_{\textup{P}}\le \mathcal C_{\textup{LS}}$. 
For $m$-strongly convex potentials, we have $\mathcal C_{\textup{SG}}^{-1} = \mathcal C_{\textup{P}}\le \mathcal C_{\textup{LS}}\le 1/m$. Results in \citep{Dalalyan14,Durmus2} 
imply that in order to get a Wasserstein distance smaller than $\varepsilon \sqrt{p/m}$, 
it suffices to perform a number of LMC iterations proportional to $(M/m)^2\varepsilon^{-2}$, up to logarithmic factors. A formal 
proof of the fact that the same result holds for the densities
satisfying the log-Sobolev inequality with constant $1/m$ (but which
are not necessarily $m$-strongly log-concave) was given in \citep{Vempala_Wibisono}. 

On the other hand, it was established by \citep{Bobkov99} that 
any log-concave distribution satisfies the Poincar\'e inequality. 
However, the Poincar\'e constant might depend  on the dimension. 
In \citep{Kannan95}, the authors conjectured that there is a universal
constant $\mathcal C_{\textup{KLS}} >0$ such that for any log-concave 
distribution $\pi$ 
on $\mathbb R^p$, 
\begin{align}\label{KLS}\tag{KLS}
    \mathcal C_{\textup{P}} \le \mathcal C_{\textup{KLS}} \|\mathbf E_\pi[\bX\bX^\top]\|_{\textup{op}}:=\mathcal C_{\textup{KLS}} \mu_{\textup{op}}(\pi). 
\end{align}
Despite important efforts made in recent years (see \citep{KLS_conjecture, cattiaux2018poincar}),
this conjecture is still unproved. Finally, in the recent paper \citep{chewi2020exponential}, Corollary 4 establishes that 
$W_2(\mu_t^{\textupsf{LD}},\pi) \le \sqrt{2\mathcal C_{\textup{P}}\chi^2(\nu_0||\pi)}\,e^{-t/\mathcal C_{\textup{P}}}$. 
While the exponential in $t$ convergence to zero is a very appealing
property of this result, it comes with two shortcomings. To the previously
mentioned difficulty of assessing the Poincar\'e constant, one has to add
the challenging problem of finding a meaningful upper bound on the 
$\chi^2$-divergence between the initial distribution and the target. 

What can we hope for in the light of the previous discussion? 
As shown in \citep[Lemma 5]{Dalalyan14}, for $f$ satisfying \SCGL{m}{M}, 
choosing $\nu_0 = \mathcal N(\bx_*,M^{-1}\bfI_p)$ yields 
$\chi^2(\nu_0\|\pi) \le (M/m)^{p/2}$. In the case $m=0$, it might be possible
to replace $m$  by $1/\mathcal C_{\textup{P}}$ in this result. If in addition, 
we admit inequality \ref{KLS}, then we get 
\begin{align}
    W_2^2(\mu_t^{\textupsf{LD}},\pi) &\le 2\mathcal C_{\textup{P}}\,(M\mathcal C_{\textup{P}})^{p/2}e^{-2t/\mathcal C_{\textup{P}}}\\
    &\le 2\mathcal C_{\textup{KLS}} \mu_{\textup{op}}(\pi)\,\big(M\mathcal C_{\textup{KLS}} \mu_{\textup{op}}(\pi)\big)^{p/2}
    e^{-2t/\mathcal C_{\textup{KLS}} \mu_{\textup{op}}(\pi)}.
\end{align}
This is, probably, the best upper bound one could hope for in the general
log-concave setting by Langevin diffusion based algorithms. We see that it 
has three drawbacks as compared to our result stated in \Cref{prop-optimal}. 
First, it requires the knowledge of a minimizer $\bx^*$. Second, it
involves the Lipschitz constant $M$ of the gradient. Third, it is heavily
based on $\mathcal{C}_{\textup{KLS}}$, which might be very large.

\section{Penalized Gradient Flow}\label{sec-C}

\subsection{Proof of \Cref{th-flow}}
We recall that for every $\gamma \in \RR $,  
 is given by $f_{\gamma}(\cdot) := f(\cdot) + \gamma \|\cdot\|^2_2/2$. We define $\bx(\gamma)$  the minimum point of 
$f_{\gamma}$. In particular, $\bx_0 = \bx_*$.  The  
triangle inequality  yields
\begin{equation}\label{eq-flfirst-tr-ineq}
    \|\bX^{\textupsf{{PGF}}}_{t} - \bx_*\|_2 \leq 	\|\bX^{\textupsf{{PGF}}}_t - \bx_{\alpha(t)}\|_2 + \|\bx_{\alpha(t)} - 
    \bx_{0}\|_2
\end{equation} 
for every $t>0$.
We will bound these two terms separately. \AssA{\sf A}{\sf q} for $\gamma = 0$ and $\tilde{\gamma} = \alpha(t)$ yields the following bound on the second term: \begin{equation}
	   \|\bx_{\alpha(t)} - \bx_{0}\|_2 \leq  {\mathsf D}\alpha(t)^{1 - {\sf q}}.
\end{equation}
To bound the first term of \eqref{eq-flfirst-tr-ineq},  we aim at obtaining a
Gronwall-type inequality for the function 
\begin{equation}
    \phi(t):= \|\bX^{\textupsf{{PGF}}}_{t} - \bx_{\alpha(t)}\|_2.
\end{equation} 
To this end, we consider an auxiliary 
stochastic process $\{\tilde{\bX}_{u}:u\ge t\}$, defined as 
a solution of the following  differential equation
\begin{equation}
    d\tilde{\bX}_{u} = -\big(\nabla f(\tilde{\bX}_{u}) + \alpha(t) \tilde{\bX}_{u}\big)du,
\end{equation}
with the starting point $\tilde{\bX}_{t} = {\bX}_{t}$.  
This is in fact the gradient flow corresponding to the strongly-convex potential
$f_{\alpha(t)}$. The triangle inequality yields
\begin{equation}
         \phi(t+\delta) 
         \leq 
         \left\|\bX^{\textupsf{{PGF}}}_{t+\delta} -\tilde{\bX}_{t+\delta}\right\|_2 +
         \left\| \tilde{\bX}_{t+\delta} - \bx_{\alpha(t)} \right\|_2 + 
         \left\|\bx_{\alpha(t)} - \bx_{\alpha(t+\delta)}  \right\|_2 .    
\end{equation}
From the linear convergence of the gradient flow of an $\alpha(t)$-strongly 
convex function, we get the following:
\begin{equation}\label{eq-contr-delta}
    \left\| \bX^{\textupsf{PGF}}_{t+\delta} - \bx_{\alpha(t)} \right\|_2  \leq 
		\exp\big(-\delta\alpha(t)\big)
	\left\| \tilde{\bX}_{t} - \bx_{\alpha(t)} \right\|_2 
		= \exp \left( -\delta \alpha(t) \right) \phi(t).
\end{equation}  
In order to bound the distance between $\bx_{\alpha(t)}$ and 
 $\bx_{\alpha(t+\delta)}$, we use again  \AssA{\sf A}{\sf q} condition, thus 
 \begin{equation}\label{eq-dist-t-t+delta}
 	\left\|\bx_{\alpha(t)} - \bx_{\alpha(t+\delta)}  \right\|_2 
 	 \leq \frac{\sf D}{\alpha^{\sf q} (t)} (\alpha(t)- \alpha(t+\delta)) 
 	 \|\bx_*\|_2.
 \end{equation}
Thus we obtain a bound  for $\phi(t+\delta)$, that depends linearly on $\phi(t)$:
\begin{align}\label{ineq-flbef-division}
	\phi(t+\delta)\leq  
	 \left\| \bX^{\textupsf{PGF}}_{t+\delta} - \tilde{\bX}_{t+\delta} \right\|_2 
	& +  e^{-\delta \alpha(t)}\phi(t) 
     + \frac{\sf D}{\alpha^{\sf q} (t)} (\alpha(t)- \alpha(t+\delta)) 
 	  \|\bx_*\|_2.
\end{align}
Let us subtract $\phi(t)$ from both sides of 
\eqref{ineq-flbef-division} and divide by $\delta$:
\begin{align}\label{eq-flminus-phi}
	\frac{\phi(t+\delta) - \phi(t)}{\delta}\leq 
		\frac{1}{\delta}\cdot \left\| \bX^{\textupsf{PGF}}_{t+\delta} - \tilde{\bX}_{t+\delta} \right\|_2 
	& + \frac{\exp\left(-\delta\alpha(t)\right) - 1}{\delta}\cdot \phi(t)\\
   & +\frac{{\sf D}(\alpha(t)- \alpha(t+\delta))}{\delta\alpha^{\sf q}(t)}
   \|\bx_*\|_2 .
\end{align}

The next lemma provides an upper bound on  
$ \left\| \bX^{\textupsf{PGF}}_{t+\delta} - \tilde{\bX}_{t+\delta} \right\|_2  $ showing that it 
is $o(\delta)$, when $\delta \rightarrow 0$. \\

\begin{lemma}\label{lem-fl-gr}
     Suppose $f$ satisfies \SCGL{m}{M} with $m=0$. Then for every  $t,\delta>0$, and for every integrable function $\alpha:[t,t+\delta]\to [0,\infty)$,
	\begin{equation}
	  \|\tilde{\bX}_{t+\delta} - \bX^{\textupsf{{PGF}}}_{t+\delta}\|_2
 	  \le \big(\phi(t) + \|\bx_{\alpha(t)}\|_2\big)
 	  \exp\bigg\{M\delta  +  \int_0^\delta \alpha(t+u)\,du \bigg\}
 	  \int_0^\delta \big|\alpha(t+s)-\alpha(t)\big|\,ds.
   \end{equation}
\end{lemma}

The proof can be found in the \Cref{sec-C-lem}.
When $\delta$ tends to $0$, according to \Cref{lem-fl-gr}, the first term of the 
right-hand side of \eqref{eq-flminus-phi} vanishes. Thus, after passing to the limit, we 
are left with the  following  Gronwall-type	inequality:
\begin{equation}\label{eq-flgron-phi}
	\phi'(t)\leq  - \alpha(t) \phi(t)  -
	\frac{{\sf D}\alpha'(t)}{ \alpha^{\sf q} (t)}\cdot 
	 \|\bx_{*}\|_2.
\end{equation}
Here we tacitly used the fact that 	 $\|\bx_{\alpha(t+\delta)}\|_2 
\leq  \|\bx_{0}\|_2$.
Recalling that
the function $\beta(t)$ is given by $\beta(t) = \int_0^t \alpha(s) ds$, one can 
rewrite  \eqref{eq-flgron-phi}  as
\begin{equation}
	\big(\phi(t) e^{\beta(t)} \big)'
	\leq -  \frac{{\sf D}\alpha'(t) e^{\beta(t)}}{\alpha^{\sf q} (t)} 
	\|\bx_*\|_2.
\end{equation}

Therefore we infer the following bound on $\phi(t)$:
\begin{align}\label{eq-flbefore-ipp}
 	\phi(t)&\leq  \phi(0) e^{-\beta(t)} -
    {\sf D}\|\bx_{*}\|_2\int_0^t \frac{\alpha'(s)}{{\alpha^{{\sf q}}(s)}}
	e^{\beta(s)-\beta(t)} ds.
\end{align}
Combining this bound with  \eqref{eq-last-term},
we obtain the inequality
\begin{align}
    \|\bX^{\textupsf{{PGF}}}_t - \bx_{0}\|_2
	& \leq   \|\bX^{\textupsf{{PGF}}}_0 - \bx(\alpha(0))\|_2 e^{-\beta(t)}  - {\sf D}\|\bx_{*}\|_2\int_0^t
	\frac{\alpha'(s)}{{\alpha^{\sf q}(s)}} e^{\beta(s)-\beta(t)} ds+
	{\sf D}\|\bx_{*}\|_2\alpha(t)^{1-{\sf q}}.
\end{align}
Since the process starts at point $0$,  $\|\bX^{\textupsf{{PGF}}}_0 - 
\bx(\alpha(0))\|_2 =  \|\bx_{\alpha(0)}\|_2$. The next lemma bounds 
$\|\bx_{\alpha(0)}\|_2$.
\begin{lemma}\label{lem-normon}
    The function $\gamma \mapsto \|\bx(\gamma)\|_2$ is a non-increasing continuous function on the interval $[0,\infty)$.
\end{lemma} 
Therefore, $\|\bx_{\alpha(0)}\|_2 \leq \|\bx_0\|_2 = \|\bx_*\|_2$, which completes the proof of \Cref{th-flow}.

\subsection{ Proof of \Cref{lem-fl-gr}}\label{sec-C-lem}

 	From the definition of $\tilde{\bX}$, we 
    can write
    \begin{equation}
        \tilde{\bX}_{t+\delta} - \bX^{\textupsf{{PGF}}}_{t+\delta} = \int_{t}^{t+\delta}
        \left(\nabla f({\bX}^{\textupsf{PGF}}_{s}) - \nabla f(\tilde{\bX}_{s}) + 
        \alpha(s){\bX}^{\textupsf{PGF}}_{s} - \alpha(t)\tilde{\bX}_{s}\right)ds.
    \end{equation}
    Therefore we have
    \begin{align}
        \|\tilde{\bX}_{t+\delta} - \bX^{\textupsf{{PGF}}}_{t+\delta}\|_{2} \leq  &
        \bigg\| \underbrace{\int_{t}^{t+\delta} \left(\nabla f({\bX}^{\textupsf{PGF}}_{s}) - 
        \nabla f(\tilde{\bX}_{s})\right)ds}_{:=T_1}\bigg\|_{2} + \bigg\| \underbrace{\int_{t}^{t+\delta} \left(\alpha(s)
        {\bX}^{\textupsf{PGF}}_{s} - \alpha(t)\tilde{\bX}_{s}\right)ds}_{:=T_2}\bigg\|_{2}. 
        \label{eq:11}
    \end{align}    Now let us analyze these two terms separately. We start with  $T_1$:
    \begin{align}
         \|T_1\|_{2} &= \left\| \int_{t}^{t+\delta} \left(\nabla f({\bX}^{\textupsf{PGF}}_{s}) - 
        \nabla f(\tilde{\bX}_{s})\right)ds\right\|_{2}\\ &\leq
         \int_{t}^{t+\delta} \left\|\nabla f({\bX}^{\textupsf{PGF}}_{s}) - 
        \nabla f(\tilde{\bX}_{s})\right\|_{2}ds \\
        & \leq M\int_{t}^{t+\delta}  \|{\bX}^{\textupsf{PGF}}_{s} - 
        \tilde{\bX}_{s}\|_{2}ds.
    \end{align}
    These are due to the Minkowskii inequality and the Lipschitz 
    continuity of the gradient.
    In order to bound the second term $T_2$, we will add and subtract 
    the term 
    $\alpha(t+s)\tilde{\bX}_{t+s}$.  Similar to the case above,
     we get the 
    following upper bound:
    \begin{align}
    	\|T_2\|_{2} &\le  \int_{t}^{t+\delta} \alpha(s) 
         \,\big\| {\bX}^{\textupsf{PGF}}_{s} - \tilde{\bX}_{s}\big\|_{2} ds 
         +  
         \int_{t}^{t+\delta} \big|\alpha(s) - \alpha(t)\big|\,\big\|\tilde{\bX}_{s}\big\|_{2}
          ds\\
          &=\int_{0}^{\delta} \alpha(t+s) 
         \,\big\| {\bX}^{\textupsf{PGF}}_{t+s} - \tilde{\bX}_{t+s}\big\|_{2} ds 
         +  
         \int_{0}^{\delta} \big|\alpha(t+s) 
         - \alpha(t)\big|\,\big\|\tilde{\bX}_{t+s} 
         \big\|_{2} ds.
    \end{align}
    Recall that 
    $\tilde{\bX}_{t+s}$ is the gradient flow  of an 
    $(m+\alpha(t))$-strongly convex  potential function. Thus, the triangle 
    inequality   yields 
    \begin{align}
    	\big\|\tilde{\bX}_{t+s}\big\|_{2}  
    	    &\le \big\| \tilde{\bX}_{t+s} - \bx(t) \big\|_2 + 
    	    \|\bx(t)\|_2 \\
    		&\leq \big\| \tilde{\bX}_{t} - \bx(t) \big\|_2
    		\exp(-ms-\alpha(t)s) +  \|\bx(t)\|_2 \\ 
    		&\leq   \big\|{\bX}^{\textupsf{PGF}}_{t} - \bx(t) \big\|_2 +  \|\bx(t)\|_2:=V_t.
    \end{align}
    Summing up, we have
    \begin{align}
         \big\|  {\bX}^{\textupsf{PGF}}_{t+\delta} - \tilde{\bX}_{t+\delta}
         \big\|_{2} &\leq
          \int_{0}^\delta \big(M +  \alpha(t+s)\big)
         \|  {\bX}^{\textupsf{PGF}}_{t+s} - \tilde{\bX}_{t+s}\|_{2}
         ds + \tilde\alpha_t(\delta)\,V_t,
	\end{align}
	where $\tilde{\alpha_t}(\delta)$ is an auxiliary function defined as
	\begin{align}
        \tilde \alpha_t(\delta) :=  \int_0^\delta
        |\alpha(t+s)-\alpha(t)|\,ds.
    \end{align}
    Now let us define  $\Phi(s) =  \|  {\bX}^{\textupsf{PGF}}_{t+s} - \tilde{\bX}_{t+s}\|_{\mathbb L_2}$. The last inequality can be rewritten as
    	\begin{equation}
    		\Phi(\delta) \leq \int_0^\delta \big(M  +  \alpha(t+s)\big)\Phi(s)\,ds
    		+\tilde\alpha_t(\delta)\, V_t.
    	\end{equation}
    The (integral form of the) Gronwall inequality  implies that
    \begin{align}
        \Phi(\delta)  &\le V_t
        \int_0^\delta 
        \tilde\alpha_t(s)\big(M  +  \alpha(t+s)\big) e^{\int_s^\delta (M  +  \alpha(t+u))\,du}\,ds + \tilde\alpha_t(\delta)V_t \\
        &= V_t \int_0^\delta 
        \tilde\alpha_t'(s)\,e^{\int_s^\delta (M  +  \alpha(t+u))\,du}\,ds\\
        &\le V_t\,\tilde\alpha_t(\delta)\,
        \exp\bigg\{M\delta  +  \int_0^\delta \alpha(t+u)\,du \bigg\}. 
    \end{align}
    This completes the proof.

\subsection{Proof of \Cref{lem-normon}}\label{sec-C-mon}

Suppose that $\gamma_1 < \gamma_2$. Let us show  that $\|\bx(\gamma_1)\|_2 > 
\|\bx(\gamma_2)\|_2$. Let us consider the function $f_{\gamma_2}$. We have that
\begin{align}
    f_{\gamma_2}(\bx(\gamma_2)) &\leq f_{\gamma_2}(\bx(\gamma_1))=f(\bx(\gamma_1)) + \gamma_2\|\bx(\gamma_1)\|_2/2.
\end{align}
The definition of $f_{\gamma_1}$  yields
\begin{align}
    f_{\gamma_2}(\bx(\gamma_2))
    &\leq f_{\gamma_1}(\bx(\gamma_1)) + (\gamma_2 - \gamma_1)\|\bx(\gamma_1)\|_2/2 \\
    &\leq f_{\gamma_1}(\bx(\gamma_2)) + (\gamma_2 - \gamma_1)\|\bx(\gamma_1)\|_2/2\\
    &\leq f_{\gamma_2}(\bx(\gamma_2)) + (\gamma_2 - \gamma_1)\big( \|\bx(\gamma_1)\|_2
    - \|\bx(\gamma_2)\|_2\big)/2.
\end{align}
Here the second passage is valid, as $\bx(\gamma_1)$ is the minimum point of 
$f_{\gamma_1}$. Since $\gamma_2 > \gamma_1$, the difference $\|\bx(\gamma_1)\|_2
- \|\bx(\gamma_2)\|_2$ is positive. Thus the monotony is proved.\\
To prove the 
continuity of the function we take a sequence $\gamma_n$ that tends to $\gamma_0$ and 
show that $\bx_{\gamma_n} \rightarrow \bx_{\gamma_0}$. Assumption \AssA{\sf D}{\sf q} 
yields
\begin{align}
      \|\bx_{\gamma_n} - \bx_{\gamma_0}\|_2 \leq
      \frac{\mathsf D}{\max(\gamma_n,\gamma_0)^{\mathsf q}} |\gamma_n - \gamma_0|
      \|\bx_*\|_2,\qquad \forall n\in \NN.
\end{align}
Since ${\sf q} < 1$, the ratio of $|\gamma_n - \gamma_0|$ and 
$\max(\gamma_n,\gamma_0)^{\mathsf q}$ tends to zero, when $n\rightarrow 0$. This concludes the proof.

\section{Examples of functions satisfying condition \protect\AssA{\sf D}{\sf q}
}

In this section we consider several functions that are convex but
not strongly convex and satisfy \AssA{\sf D}{\sf q} condition presented
in \Cref{sec:PGF}. 

\subsection{Locally strongly convex functions}
We prove that locally strongly convex functions satisfy \AssA{\sf D}{0}.
Recalling \Cref{lem-normon} we get that $\|\bx_\gamma\|_2 \leq \|\bx_* \|_2$. Thus the
we can consider the function only on $\mathcal{B}(0,\|x_*\|_2)$. Since $f$ is locally 
strongly convex, there exists $m_*$ such that it is $m_*$-strongly convex in the ball
$\mathcal{B}(0,\|x_*\|_2)$. The latter means, that $f_{\tilde{\gamma}}$ is
$(m_* + \tilde{\gamma})$-strongly convex. Therefore \citep{Nest}[Theorem 2.1.9] yields
the following:
\begin{equation}
    \|\bx_{\gamma} - \bx_{\tilde\gamma}\|_2 \leq \frac{1}{m_* + {\tilde\gamma}}
    \|\nabla f_{\tilde\gamma}(\bx_{\gamma})  - \nabla 
    f_{\tilde\gamma}(\bx_{\tilde\gamma}) \|_2 .
\end{equation}
Using the optimality condition for differentiable functions one gets 
$\nabla f_{\tilde\gamma}(\bx_{\gamma}) = ({\tilde\gamma} - \gamma) \bx_\gamma$ for all
$\gamma \geq 0$. Therefore, for every $0 \leq \gamma < {\tilde{\gamma}}$, we obtain
\begin{equation}
    \|\bx_{\gamma} - \bx_{\tilde\gamma}\|_2 \leq \frac{1}{m_* + {\tilde\gamma}}
    \| ({\tilde\gamma} - \gamma) \bx_\gamma \|_2
    \leq \frac{{\tilde\gamma} - \gamma}{m_*} \|\bx_\gamma \|_2. 
\end{equation}
The latter is true due to \Cref{lem-normon}. Thus $f$ satisfies \AssA{1/m_*}{0}.

\subsection{Cubic function $f(\bx) = \|\bx - \bx_*\|_2^3$}

In this section we show that the cubic function satisfies 
\AssA{\nicefrac{1}{\sqrt{3\|\bx_*\|_2}}}{\nicefrac{1}{2}}. 
It is straightforward to verify that the function $f$ is convex. $f_\gamma$ is 
strongly convex and the optimality condition for $\bx_\gamma$ yields the following 
equality:
\begin{align}\label{eq-pengr0}
    \nabla f(\bx_\gamma) + \gamma \bx_\gamma = 3\|\bx_\gamma - \bx_*\|_2(\bx_\gamma - 
    \bx_*) + \gamma \bx_\gamma = 0.
\end{align}
In the case when $\bx_* = 0$, the penalized minimum point $\bx_\gamma$ equals $0$, for
every $\gamma$, thus we suppose in the following that $\bx_* \neq 0$. 
Since the norm is scalar, \eqref{eq-pengr0} yields that the vectors $\bx_\gamma - 
\bx_*$ and $\bx_\gamma$ are co-linear. Therefore there exists a real number 
$\lambda_\gamma$ such that $\bx_\gamma = \lambda_\gamma \bx_*$. \Cref{lem-normon} 
implies that $|\lambda_\gamma| \leq 1$,  thus the following quadratic equality is 
true:
\begin{align}\label{eq-grlambda}
    - 3\|\bx_*\|_2(\lambda_\gamma - 1)^2 \bx_* + \gamma \lambda_\gamma \bx_*  = 0.
\end{align}
As said in the beginning, $\bx_* \neq 0$, therefore it its coefficient that is equal 
to zero. Solving the quadratic equation with respect to $\lambda_\gamma$, we get the 
following formula:
\begin{align}
    \lambda_\gamma = 1 - \frac{\gamma}{\gamma/2 +
    \sqrt{3\gamma\|\bx_*\|_2 + \gamma^2/4}}.
\end{align}
According to \Cref{lem-normon}, for every $\tilde{\gamma} > \gamma$, we have 
$|\lambda_\gamma| > |\lambda_{\tilde\gamma}|$. On the other hand, from 
\eqref{eq-grlambda} one deduces that $\lambda_\gamma > 0$, for every $\gamma > 0$. 
Thus, inserting the found value for
$\lambda_\gamma$, we obtain the following 
inequality:
\begin{align}
    \|\bx_\gamma - \bx_{\tilde\gamma}\|_2 &= \|\bx_*\|_2 
    \left(\frac{\tilde\gamma}{{\tilde\gamma}/2 +
    \sqrt{3{\tilde\gamma}\|\bx_*\|_2 + {\tilde\gamma}^2/4}} - \frac{\gamma}{\gamma/2 +
    \sqrt{3\gamma\|\bx_*\|_2 + \gamma^2/4}}\right) \\
    &\leq \frac{ ({\tilde\gamma} - \gamma)  \|\bx_*\|_2 }{{\tilde\gamma}/2 +
    \sqrt{3{\tilde\gamma}\|\bx_*\|_2 + {\tilde\gamma}^2/4}}
    \leq \frac{ ( {\tilde\gamma} - \gamma)  \|\bx_*\|_2 }
    {\sqrt{3{\tilde\gamma}\|\bx_*\|_2 }}.
\end{align}
Therefore $f$ satisfies \AssA{\nicefrac{1}{\sqrt{3\|\bx_*\|_2}}}{\nicefrac{1}{2}}.

\subsection{Power function $f(\bx) = \|\bx-\bx_*\|^a_2$}

For $a\ge 2$, we consider the function $f(\bx) = \|\bx-\bx_*\|^a_2$.  We show here 
that $f$ satisfies \AssA{(\nicefrac{1}{a 
\|\bx_*\|_2^{a-2}})^{1/(a-1)}}{\nicefrac{(a-2)}{(a-1)}}.
Since $f_\gamma$ is a differentiable strongly-convex function, we get the following equation for $\bx_\gamma$:
\begin{align}\label{eq-powgrad}
    a\|\bx_\gamma - \bx_*\|_2^{a-2} (\bx_\gamma - \bx_*) + \gamma\bx_\gamma = 0.
\end{align}
Similar to the previous case, we notice that $\bx_\gamma - \bx_*$ and $\bx_\gamma$ are
co-linear. Thus, there exists $\lambda_\gamma$ such that $\bx_\gamma = 
(1 - \lambda_\gamma)\bx_*$. Since $\bx_*$ is assumed to be non-zero, in order to calculate 
$\bx_\gamma$, one needs to solve the following equation:
\begin{align}\label{eq-powlamb}
    |\lambda_\gamma |^{a-2}\lambda_\gamma = \frac{\gamma (1 - \lambda_\gamma)}{a \|\bx_*\|^{a-2}}.
\end{align}
Thus the  $p$-dimensional equation \eqref{eq-powgrad} reduces to equation 
\eqref{eq-powlamb} involving a one-dimensional unknown. \Cref{lem-normon} yields 
$\lambda_{\tilde\gamma} > \lambda_\gamma>0 $ for every
$\tilde\gamma > \gamma \geq 0$. In addition, from \eqref{eq-powlamb}, we have that 
$\lambda_\gamma\le 1$ for every $\gamma > 0$. It is straightforward to verify 
that for every $\gamma \geq 0$, \eqref{eq-powlamb} has exactly one solution satisfying 
these conditions.

\begin{lemma}
Let $\alpha\ge 1$. If  $(\lambda_s:s\in(0,1))$ satisfies $\lambda_s^\alpha = s(1-\lambda_s)$ for every $s\in(0,1)$, 
then 
\begin{align}
    |\lambda_s - \lambda_{s'}| \le \frac{|s-s'|}{(s\vee s')^{(\alpha-1)/\alpha}},\qquad \forall s', s\in (0, 1).
\end{align}
\end{lemma}
\begin{proof}
Without loss of generality, we assume that $s'\le s$. Computing
the derivative of both sides of the identity $\lambda_s^\alpha = s(1-\lambda_s)$, we get 
\begin{align}
    \lambda'_s = \frac{1-\lambda_s}{\alpha \lambda_s^{\alpha-1}+s}
    \ge 0.
\end{align}
This implies that $\lambda_{s'}\le \lambda_s$. In addition,
\begin{align}
    \lambda_{s}-\lambda_{s'} &\le \frac{\lambda_{s}^\alpha-\lambda_{s'}^\alpha}
    {\lambda_{s}^{\alpha-1}}\\ 
    & = \frac{s(1-\lambda_s)-s'(1-\lambda_{s'})}{\lambda_s^{\alpha-1}}\\
    & = \frac{(s-s') (1-\lambda_{s'})}{\lambda_{s}^{\alpha-1}} - \frac{s(\lambda_{s}-\lambda_{s'})}{\lambda_{s}^{\alpha-1}}.
\end{align}
Rearranging the terms, we arrive at
\begin{align}
    \lambda_{s}-\lambda_{s'} 
    &\le \frac{(s-s') (1-\lambda_{s'})}{\lambda_{s}^{\alpha-1}
    \big(1 + \frac{s}{\lambda_{s}^{\alpha-1}}\big)}\\
    & = \frac{(s-s') (1-\lambda_{s'})}{\lambda_{s}^{\alpha-1}
    + s}
\end{align}
In the last fraction, the numerator is bounded by $s-s'$, while
the denominator satisfies
\begin{align}
    \lambda_{s}^{\alpha-1}
    + s & = (s(1-\lambda_{s}))^{(\alpha-1)/\alpha} + s\\
    &\ge (s(1-s^{1/\alpha}))^{(\alpha-1)/\alpha} + s\\
    &\ge s^{(\alpha-1)/\alpha}(1-s^{1/\alpha}) + s = 
    s^{(\alpha-1)/\alpha}. 
\end{align}
This completes the proof of the lemma. 
\end{proof}
Applying this lemma to \eqref{eq-powlamb}, we get
\begin{align}
    \lambda_{\tilde\gamma} -\lambda_\gamma \le 
    \frac{\tilde\gamma - \gamma }{a\|\bx_*\|_2^{a-2}(\tilde\gamma /a\|\bx_*\|_2^{a-2})^{(a-2)/(a-1)}}
    =\frac{\tilde\gamma - \gamma }{a^{1/(a-1)}\|\bx_*\|_2^{(a-2)/(a-1)}\tilde\gamma ^{(a-2)/(a-1)}},
\end{align}
for all $\gamma,\tilde\gamma$ satisfying 
$0 \le \gamma \le \tilde\gamma \le a \|\bx_*\|_2^{a-2}$. 
In conclusion, we get
\begin{align}
    \|\bx_{\tilde\gamma} - \bx_{\gamma}\|_2 
    &\le \frac{\tilde\gamma - \gamma }{a^{1/(a-1)}\|\bx_*\|_2^{(a-2)/(a-1)}\tilde\gamma ^{(a-2)/(a-1)}}
    \|\bx_*\|_2\\
    &\leq \frac{\tilde\gamma - \gamma }{\tilde\gamma ^{(a-2)/(a-1)}}
    (\|\bx_*\|_2/a)^{1/(a-1)} .
\end{align}
This concludes the proof.

\end{document}